\documentclass{article}
\usepackage{arxiv}

\usepackage{graphicx} 
\usepackage{amssymb,amsthm,amsmath,amsfonts,mathtools,latexsym}

\theoremstyle{plain}
\newtheorem{theorem}{Theorem}[section]

\newtheorem{lemma}[theorem]{Lemma}

\theoremstyle{definition}
\newtheorem{definition}[theorem]{Definition}
\newtheorem{example}[theorem]{Example}

\theoremstyle{remark}
\newtheorem{remark}[theorem]{Remark}

\usepackage{dirtytalk} 
\usepackage{url} 
\usepackage{wasysym}
\DeclareFontFamily{U}{musix}{}
\DeclareFontShape{U}{musix}{m}{n}{
	<-12>   musix11
	<12-15> musix13
	<15-18> musix16
	<18-23> musix20
	<23->   musix29
}{}
\newcommand*\musix{\usefont{U}{musix}{m}{n}\selectfont}
\DeclareTextFontCommand{\textmusix}{\musix}
\newcommand*\doubleflat{\raisebox{.6ex}{\textmusix{3}}}
\newcommand*\doublesharp{\raisebox{.6ex}{\textmusix{5}}}
\usepackage{relsize}

\usepackage{tikz} 
\usetikzlibrary{arrows}
\usepackage{xcolor}

\usepackage{authblk}

\author[1]{Caleb Scott Alons\thanks{Email: ca6w6@umsystem.edu}}
\affil[1]{Department of Mathematics, University of Missouri - Columbia}

\begin{document}

\title{A General Approach to Enharmonicism}

\maketitle

\begin{abstract}
	We propose a general approach to enharmonicism within syntactic music theory; that is, we formalize several abstract structures of enharmonicism without any acoustic considerations. The structures of concern to this paper are constructed from diatonics, chromatics, and enharmonics, and we seek to present a general enharmonic theory that builds upon Hook's enharmonic theory published in \say{Enharmonic Systems: A Theory of Key Signatures, Enharmonic Equivalence and Diatonicism} \cite{Hook}. We also propose that the two criteria of reflection and maximal evenness can be used to determine whether an arbitrary enharmonic system is musically practical. Moreover, we argue that preferential treatment of the standard enharmonic system (SES) can be motivated from an abstract mathematical perspective divorced from pitch consideration and conventional acoustic constructions of SES.
\end{abstract}

\section{Context and Motivation}
\label{SEC.Context and Motivation}

\subsection{Historical Considerations}
\label{SSEC.Historical Considerations}

The study of enharmonicism has been historically treated as a primarily acoustic problem. Before the development of equal temperament, accidentals were introduced because it was not possible to provide a complete set of justly intoned concords for each of the $7$ diatones. During the Renaissance, various attempts were made to explore enharmonicism more generally. Nicola Vicentino constructed a keyboard called the Archicembalo, which had $36$ physical keys per octave to realize a $31$-tone pitch system \cite{Vicentino}. Guillaume Costeley's 1558 composition, \say{Seigneur Dieu ta pitié s'estende dessus moy}, required a special $19$-tone tuning, and Costeley explicitly intended for his $19$-tone system to form a cyclic equal division of the octave---an approach differing from Vicentino in the sense that Vicentino's pitch divisions were expressly nonuniform within the octave \cite{Costeley,Wibberley}. In his 1577 \textit{De musica libri septem}, Francisco de Salinas mathematically described a $1/3$-comma meantone system with $19$ pitches per octave, whose intervallic structure closely approximates $19$-TET, while also describing a $24$-pitch system based on just intonation \cite{Salinas}. Later in 1627, Mersenne proposed mathematical and acoustic bases of concepts like consonance, interval division, and temperament, where his text gives a theoretical framework for tuning expanded pitch systems \cite{Mersenne}.

Note that, for all the above, standard equal temperament tuning had not been established yet, hence why they did not consider transposability, and why accidentals played intentional acoustic roles within enharmonicism. However, once standard equal temperament became widely accepted, accidentals played a more syntactic role with respect to tonal transposability. At the start of the 20th century, however, composers such as Harry Partch and Ben Johnston revived the expansion of standard enharmonicism through justly intoned intervallic relationships, deriving expanded pitch vocabularies from rational frequency ratios rather than from the equal subdivision of the octave \cite{Partch,JohnstonProportionality,johnston_notation_2006}. There even exist acoustic motivations for using extended sets of diatones in recent history; Cubarsi, for example, acoustically motivates the use of extended sets of diatones by showing how successive convergents and semiconvergents of the Pythagorean generator yield increasingly fine microtonal subdivisions, allowing tones to be selected according to their proximity to desired frequency ratios \cite{Cubarsi}.

We observe that it took centuries for the standard $12$-TET system to be accepted as the more-or-less universal enharmonic system, arising primarily from acoustic considerations. In this paper, we conjecture that there is also a purely mathematical reason supporting the preference of standard $12$-tone enharmonicism. To accomplish this goal, we develop a general approach to enharmonicism divorced from pitch considerations, as motivated by Hook's acknowledgment of the merit in such an approach to his enharmonic theory \cite{Hook}. However, whereas Hook constructs a syntactic enharmonic theory from the $7$ standard diatones, we construct a general syntactic theory from an arbitrary number of diatones. By starting with arbitrary diatonic cardinality, we avoid presupposing any standard acoustic preferences.

\subsection{Acoustics and Syntactics}
\label{SSEC.Acoustics and Syntactics}

We make some necessary remarks to motivate our differentiation of syntactics and acoustics, where \textit{syntactics} refers to the abstract formal structures and syntax of music theory having no inherent reference to pitch and \textit{acoustics} refers to the relationships and interactions between sound waves derived from the physical properties of sound waves themselves. Although there is little interest in comprehensively classifying every concept within music theory as either syntactic or acoustic, the distinction between syntactic and acoustic music theory is recognized both implicitly \cite{Clough_Douthett} and explicitly \cite{Hook}.

To the best of the author's knowledge, a historical tendency of music theory has been to derive syntactics from acoustics: through a variety of methods and experiments, musicians and composers identify useful acoustic phenomena---such as justly intoned intervals, concords, the overtone series, or particular temperaments---and subsequently seek a logical and easily transmittable syntax through which these phenomena and their musical applications can be notated and reproduced \cite{Tinctoris,Zarlino,Rameau}. This longstanding acoustic-to-syntactic relationship has, with time, encouraged the conflation of acoustics and syntactics, since syntactic systems serve the practical necessities of preserving and reproducing musical works \cite{Guido}.

However, within this history were also those who recognized how theoretical treatments of syntactics could become increasingly detached from acoustics \cite{Fetis}, and in the 21st century, we identify Hook's enharmonic theory to be a prime example of treating syntactics as an abstract study in its own right. We suggest thinking of music as the marriage of mathematics (syntactics) and sound (acoustics), and this paper shows that a general approach to enharmonicism may reveal deeper reasons for the predominant use of standard enharmonicism.

\section{Diatonic Theory}
\label{SEC.Diatonic Theory}

\subsection{Diatonic Sets}
\label{SSEC.Diatonic Sets}

Rather than beginning with the conventional assumption of the $7$ standard diatones C, D, E, F, G, A, and B, we define an arbitrary number of diatones and collect them in a set.

\begin{definition}
\label{DEF.Diatonic Set}
	A finite set $D=\{d_0,\dots,d_{n-1}\}$ with indices in $\mathbb{Z}_n$ is called a \textit{diatonic set}, the elements of a diatonic set are called \textit{diatonics} (or \textit{diatones}), and the number $n$ of diatonics is called the \textit{diatonic cardinality}.
\end{definition}

\begin{definition}
	The diatonic set of diatones $d_0=\textrm{C}$, $d_1=\textrm{D}$, $d_2=\textrm{E}$, $d_3=\textrm{F}$, $d_4=\textrm{G}$, $d_5=\textrm{A}$, and $d_6=\textrm{B}$ is the \textit{standard diatonic set}, denoted $\mathcal{D}_7$.
\end{definition}

\subsection{Diatonic Transpositions and Prime Diatonic Orderings}
\label{SSEC.Diatonic Transpositions and Prime Diatonic Orderings}

We formalize the notion of transposing between the diatones of a diatonic set.

\begin{definition}
	Let $D$ be a diatonic set and $x\in\mathbb{Z}$; a function $t_x:D\to D$ defined by the equation $t_x(d_i)=d_{i+x}$ is called the \textit{diatonic transposition} by $x$. If the diatonic cardinality of $D$ is $n$, then the index $i+x$ of the diatonic $d_{i+x}\in D$ is understood to be $i+x\pmod{n}$.
\end{definition}

\begin{example}
	Let $D_{19}=\{d_0,\dots,d_{18}\}$ be a diatonic set. Then the diatonic transposition by $10$ of $d_{14}$ is $t_{10}(d_{14})=d_{14+10}=d_5$ because $24\equiv 5\pmod{19}$, and the diatonic transposition by $-7$ of $d_4$ is $t_{-7}(d_4)=d_{4+(-7)}=d_{16}$ because $-3\equiv 16\pmod{19}$.
\end{example}

\begin{definition}
	Given a diatonic set $D$, we call the diatonic transposition by $1$ the \textit{diatonic step transposition}. Given diatones $d_i,d_j\in D$, we say that $d_j$ is $y$ diatonic steps \textit{above} $d_i$ (or, equivalently, that $d_i$ is $y$ diatonic steps \textit{below} $d_j$) if and only if $$\underbrace{(t_1\circ\cdots\circ t_1)}_y(d_i)=t_y(d_i)=d_j.$$
\end{definition}

We formalize the notion of diatonic intervals.

\begin{definition}
\label{DEF.Interval}
	Let $D$ be a diatonic set, $p$ be a positive integer, and $d_j\in D$ be $p-1$ diatonic steps above $d_i\in D$. A \textit{diatonic interval} is an interval from $d_i$ to $d_j$ in $D$, denoted $[d_i,d_j]_D$; the \textit{length} of $[d_i,d_j]_D$ is $p$, and we call $[d_i,d_j]_D$ a \textit{$p$th}. In the special case $p=1$, we call the interval $[d_i,d_i]_D$ a \textit{unison} of trivial length $0$.
\end{definition}

\begin{example}
	Since A is $2$ diatonic steps above F in $\mathcal{D}_7$, the interval length from F to A is $3$, and the interval $[\textrm{F},\textrm{A}]_{\mathcal{D}_7}$ is called a third. Since Definition \ref{DEF.Interval} accommodates the cyclic nature of diatonic sets, then C is $6$ diatonic steps above D, so the interval length from D to C is $7$, and the interval $[\textrm{D},\textrm{C}]_{\mathcal{D}_7}$ is called a seventh.
\end{example}

We turn to the concept of generating a diatonic ordering from a fixed interval length. Cyclic orderings of diatonic sets are important structures of diatonicism. Two conventionally assumed orderings of $\mathcal{D}_7$ are already familiar: the \say{scalar} ordering C, D, E, F, G, A, B (interval length $2$) and the \say{circle of fifths} ordering F, C, G, D, A, E, B (interval length $5$) \cite{Hook}. We depart from $\mathcal{D}_7$ and formalize the notion of ordering a diatonic set using an arbitrary interval length.

\begin{definition}
\label{DEF.PDO}
	Let $D=\{d_0,\dots,d_{n-1}\}$ be a diatonic set with $n>1$. Given an interval length $p$, such that $p-1$ is coprime with $n$, let $\delta_p:D\to D$ be a function defined by the equation $\delta_p(d_x)=d_{(p-1)x}$. Elements in the range of $\delta_p$ are denoted $d_i'$, where $d_i'=\delta_p(d_i)$, and are called \textit{prime diatonics}. An ordering $d_0',d_1',\dots,d_{n-1}'$ of $D$ is called a \textit{prime diatonic ordering} (\textit{PDO}), and if $p$ is known, we call the resulting PDO a \textit{diatonic circle of $p$ths}.
\end{definition}

Our use of $\delta_p(d_x)=d_{(p-1)x}$ rather than $\delta_p(d_x)=d_{px}$ to produce a diatonic circle of $p$ths accommodates the awkwardness of traditional diatonic interval numbering (where a \say{fifth} is $4$ diatonic steps, for example), so that the familiar \say{circle of fifths} terminology may be smoothly extended to \say{circle of $p$ths.} Furthermore, the requirement that $p-1$ be coprime with $n$ has been well known to music theorists (albeit articulated primarily in acoustic contexts) \cite{Gamer}, so we find Definition \ref{DEF.PDO} to be a natural way to generalize the notion of cyclic orderings of diatones.

The function $\delta_p$ is invertible, and $\delta_p^{-1}:D\to D$ is defined by the equation $\delta_p^{-1}(d_x)=d_{(p-1)^{-1}x}$, where $(p-1)^{-1}$ denotes the multiplicative inverse of $p-1$ in $\mathbb{Z}_n$. By Definition \ref{DEF.PDO}, the interval length $p$ is chosen such that $p-1$ is coprime with $n$; furthermore, a corollary of Bézout's identity states that if $a\in\mathbb{Z}_n$ and $\gcd{(a,n)}=1$, then there exists a multiplicative inverse of $a$ in $\mathbb{Z}_n$ \cite{Bezout}. Thus, $p-1$ is guaranteed a multiplicative inverse $(p-1)^{-1}$ in $\mathbb{Z}_n$. We computationally verify that $\delta_p^{-1}$ is the inverse function of $\delta_p$ as follows: \begin{align*} \delta_p(\delta_p^{-1}(d_x)) &= \delta_p(d_{(p-1)^{-1}x})=d_{(p-1)(p-1)^{-1}x}=d_x; \\ \delta_p^{-1}(\delta_p(d_x)) &= \delta_p^{-1}(d_{(p-1)x})=d_{(p-1)^{-1}(p-1)x}=d_x. \end{align*} These functions are computationally useful because they can be condensed using the prime notation previously introduced: $d_x'=d_{(p-1)x}$; $d_x=d_{(p-1)^{-1}x}'$. We briefly remark that we may express diatonic transpositions with respect to prime diatonics using these equations, where $t_x(d_i')=t_x(d_{(p-1)i})=d_{(p-1)i+x}=d_{(p-1)^{-1}[(p-1)i+x]}'=d_{i+(p-1)^{-1}x}'$.

\begin{example}
	We identify the conventionally assumed diatonic circle of fifths of $\mathcal{D}_7$ as a PDO. For a diatonic circle of fifths ($p=5$), we compute the prime diatonics using $d_x'=d_{4x}$. We show that \begin{align*} d_0' &= d_{4(0)}=d_0=\textrm{C}, & d_3' &= d_{4(3)}=d_5=\textrm{A}, & d_5' &= d_{4(5)}=d_6=\textrm{B}, \\ d_1' &= d_{4(1)}=d_4=\textrm{G}, & d_4' &= d_{4(4)}=d_2=\textrm{E}, & d_6' &= d_{4(6)}=d_3=\textrm{F}, \\ d_2' &= d_{4(2)}=d_1=\textrm{D}, \end{align*} so the diatonic circle of fifths of $\mathcal{D}_7$ is C, G, D, A, E, B, F.
\end{example}

\begin{definition}
	The diatonic circle of fifths of $\mathcal{D}_7$ is the \textit{standard diatonic circle of fifths} (\textit{standard PDO}).
\end{definition}

We briefly note that, although F, C, G, D, A, E, B and C, G, D, A, E, B, F are technically distinct orderings of $\mathcal{D}_7$, their intervals (fifths) are the same, and either can be expressed as a diatonic transposition of the other; that is, $t_4(\textrm{F}),t_4(\textrm{C}),\dots,t_4(\textrm{B})$ is just $\textrm{C},\textrm{G},\dots,\textrm{F}$. Generally speaking, any diatonic circle of $p$ths may be transposed to any desired \say{starting} diatonic; this idea receives formal consideration in Definition \ref{DEF.Diatonic Space}.

\begin{example}
	Let $D_{10}=\{d_0,\dots,d_9\}$ be a diatonic set. Since $n=10$, the only intervals that produce PDOs of $D_{10}$ are seconds, fourths, eighths, and tenths, since $2-1=1$, $4-1=3$, $8-1=7$, and $10-1=9$ are all coprime with $10$. On the other hand, thirds, fifths, sixths, sevenths, and ninths are intervals that do not produce PDOs of $D_{10}$, since $3-1=2$, $5-1=4$, $6-1=5$, $7-1=6$, and $9-1=8$ are not coprime with $10$. For a diatonic circle of seconds ($p=2$), we compute the prime diatonics using $d_x'=d_x$, and the resulting PDO is $d_0,d_1,d_2,\dots,d_9$. For a diatonic circle of fourths ($p=4$), we compute the prime diatonics using $d_x'=d_{3x}$, and the resulting PDO is $d_0,d_3,d_6,d_9,d_2,d_5,d_8,d_1,d_4,d_7$.
\end{example}

We briefly comment on PDOs produced by interval lengths that exceed the diatonic cardinality $n$ of a diatonic set. Although Definition \ref{DEF.PDO} computationally accommodates this, such PDOs are redundant, since the resulting PDO can be produced by an interval length that does not exceed $n$.

\begin{example}
	The diatonic circle of fourteenths ($p=14$) of $D_{10}$ is $d_0,d_3,d_6,\dots,d_7$, which we already showed is the diatonic circle of fourths, and $4\equiv 14\pmod{10}$.
\end{example}

In general, if an interval length $p$ exceeds the diatonic cardinality $n$, then the resulting PDO is the same as the PDO produced from the interval length $p\pmod{n}$.

\begin{definition}
\label{DEF.SDO}
	A diatonic circle of seconds is a \textit{scalar diatonic ordering} (\textit{SDO}). The SDO of $\mathcal{D}_7$ is the \textit{standard SDO}.
\end{definition}

Definition \ref{DEF.SDO} generalizes the conventional scalar diatonic ordering C, D, E, F, G, A, B of $\mathcal{D}_7$ and is motivated by the fact that every diatonic set containing at least $2$ diatonics has a diatonic circle of seconds ($p=2$, and $p-1=1$ is trivially coprime with all integers greater than $1$) produced by $d_x'=d_x$.

\subsection{Diatonic Systems}
\label{SSEC.Diatonic Systems}

As previously noted in Section \ref{SSEC.Diatonic Transpositions and Prime Diatonic Orderings}, given a PDO of a diatonic set, we may transpose each diatonic by a fixed amount of steps to obtain a permutation of the original PDO. Since we want to recognize each of these \say{rotations} (which we formalize in Definition \ref{DEF.Rotation}) as distinct (though related) objects, we formalize the notion of a diatonic space, which orders a diatonic set with a circle of $p$ths and equips it with a diatonic transposition.

\begin{definition}
\label{DEF.Diatonic Space}
	A \textit{diatonic space} $\mathfrak{d}_{p,\alpha}^D$ is a diatonic set $D$ ordered by a circle of $p$ths and equipped with a diatonic transposition $t_\alpha$.
\end{definition}

\begin{definition}
\label{DEF.Rotation}
	A diatonic space $\mathfrak{d}_{p,\alpha}^D$ is a \textit{rotation} of a diatonic space $\mathfrak{d}_{q,\beta}^{D'}$ if and only if $D=D'$ and $p=q$; the diatonic space $\mathfrak{d}_{p,\alpha}^D$ is the \textit{trivial rotation} of $\mathfrak{d}_{p,\alpha}^D$. A set containing a diatonic space $\mathfrak{d}_{p,\alpha}^D$ and all of its rotations is called a \textit{rotation class}, denoted $\mathfrak{d}_p^D$. We call $\mathfrak{d}_p^D$ a \textit{class of $p$ths}.
\end{definition}

\begin{definition}
\label{DEF.Diatonic System}
	The \textit{diatonic system} of a diatonic set $D=\{d_0,\dots,d_{n-1}\}$ is the collection of rotation classes $\mathfrak{d}^D=\{\mathfrak{d}_p^D\,|\,\gcd{(p-1,n)}=1\}$.
\end{definition}

\begin{definition}
\label{DEF.SDS}
	The diatonic system of $\mathcal{D}_7$ is the \textit{standard diatonic system} (\textit{SDS}). The rotation classes of $\mathfrak{d}^{\mathcal{D}_7}$ are the \textit{standard rotation classes}.
\end{definition}

\begin{example}
	Recall our previous comments about distinguishing C, G, D, A, E, B, F and F, C, G, D, A, E, B. Definition \ref{DEF.Diatonic Space} tells us that we have two distinct diatonic spaces, $\mathfrak{d}_{5,0}^{\mathcal{D}_7}=\{\textrm{C},\textrm{G},\textrm{D},\textrm{A},\textrm{E},\textrm{B},\textrm{F}\}$ and $\mathfrak{d}_{5,3}^{\mathcal{D}_7}=\{\textrm{F},\textrm{C},\textrm{G},\textrm{D},\textrm{A},\textrm{E},\textrm{B}\}$, and Definition \ref{DEF.Rotation} tells us that $\mathfrak{d}_{5,0}^{\mathcal{D}_7}$ and $\mathfrak{d}_{5,3}^{\mathcal{D}_7}$ are rotations of each other. Figure \ref{FIG.StandardClassOfFifthsRotations} illustrates $\mathfrak{d}_{5,0}^{\mathcal{D}_7}$ and $\mathfrak{d}_{5,3}^{\mathcal{D}_7}$ and visually motivates our \say{rotation} terminology used to describe diatonic spaces in the same rotation class.
	
	\begin{figure}[htp]
		\begin{center}
			\includegraphics[scale=0.6]{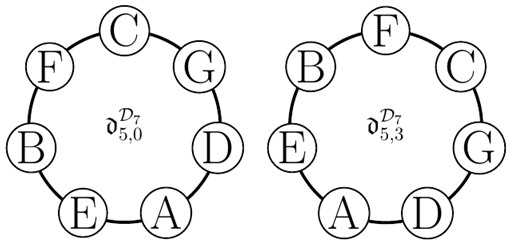}
			\caption{Some rotations in the standard class of fifths.}
			\label{FIG.StandardClassOfFifthsRotations}
		\end{center}
	\end{figure}
	
	Since the diatonic cardinality of $\mathcal{D}_7$ is $7$, the rotation class $\mathfrak{d}_5^{\mathcal{D}_7}$ contains the rotations $\mathfrak{d}_{5,0}^{\mathcal{D}_7},\mathfrak{d}_{5,1}^{\mathcal{D}_7},\dots,\mathfrak{d}_{5,6}^{\mathcal{D}_7}$. We note further that the diatonic cardinality of $\mathcal{D}_7$ is prime itself, so every interval length $p\le 7$ produces a $p-1$ coprime with $7$. The diatonic space $\mathfrak{d}_{2,0}^{\mathcal{D}_7}=\{\textrm{C},\textrm{D},\textrm{E},\textrm{F},\textrm{G},\textrm{A},\textrm{B}\}$ and its $7$ rotations form a class of seconds, the diatonic space $\mathfrak{d}_{3,0}^{\mathcal{D}_7}=\{\textrm{C},\textrm{E},\textrm{G},\textrm{B},\textrm{D},\textrm{F},\textrm{A}\}$ and its $7$ rotations form a class of thirds, and so forth, and the collection of rotation classes $\mathfrak{d}^{\mathcal{D}_7}=\{\mathfrak{d}_2^{\mathcal{D}_7},\mathfrak{d}_3^{\mathcal{D}_7},\dots,\mathfrak{d}_7^{\mathcal{D}_7}\}$ is SDS. Figure \ref{FIG.RepresentativeStandardRotationClasses} shows representative diatonic spaces from each of the six standard rotation classes.
	
	\begin{figure}[htp]
		\begin{center}
			\includegraphics[scale=0.75]{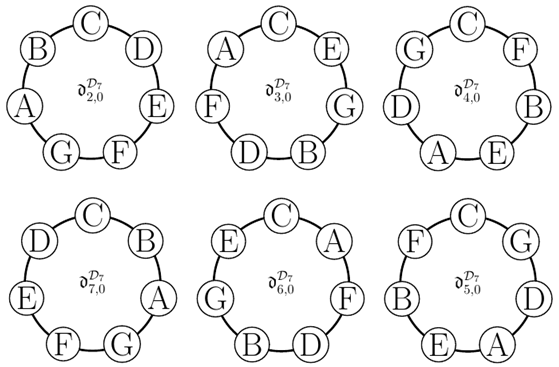}
			\caption{Diatonic spaces representing each of the six standard rotation classes.}
			\label{FIG.RepresentativeStandardRotationClasses}
		\end{center}
	\end{figure}
\end{example}

Note how $\mathfrak{d}_{2,0}^{\mathcal{D}_7}$ and $\mathfrak{d}_{7,0}^{\mathcal{D}_7}$ are mirror reflections of each other. The standard class of seconds corresponds to the unit $1$ in $\mathbb{Z}_7$ and the standard class of sevenths corresponds to the unit $6$ in $\mathbb{Z}_7$, and it is no coincidence that $1$ and $6$ are each other's negatives in $\mathbb{Z}_7$. Similar observations are made for the other standard rotation classes that are mirror reflections of each other. In short, diatonicism is the application of $\mathbb{Z}_n$ and the multiplicative group of units $U(\mathbb{Z}_n)$, where diatonic spaces and rotation classes are isomorphic to $\mathbb{Z}_n$ and diatonic systems are isomorphic to $U(\mathbb{Z}_n)$.
	
The visually intuitive notion that a diatonic system can contain \say{reflected} pairs of rotation classes proves to be significant for our later development of enharmonic systems, so we conclude our discussion of diatonicism by formalizing some terminology.

\begin{definition}
\label{DEF.Diatonic Reflection}
	Let $\mathfrak{d}^D$ be a diatonic system of a diatonic set $D$, where $|D|=n$. A rotation class $\mathfrak{d}_p^D\in\mathfrak{d}^D$ is a \textit{diatonic reflection} of a rotation class $\mathfrak{d}_q^D\in\mathfrak{d}^D$ if and only if $p+q-2=n$; that is, if $p-1$ and $q-1$ are additive inverses in $\mathbb{Z}_n$.
\end{definition}

\begin{definition}
\label{DEF.Reflectable}
	If a rotation class $\mathfrak{d}_p^D$ has a reflection in $\mathfrak{d}^D$, we call $\mathfrak{d}_p^D$ \textit{reflectable}. The relation $\sim_r$ on $\mathfrak{d}^D$, defined by $\mathfrak{d}_p^D\sim_r\mathfrak{d}_q^D$ if and only if $\mathfrak{d}_q^D$ is a reflection of $\mathfrak{d}_p^D$, is called the \textit{reflection relation}.
\end{definition}

\begin{definition}
\label{DEF.Reflection Property}
	A diatonic system $\mathfrak{d}^D$ is said to have the \textit{reflection property} (\textit{RP}) if and only if $\mathfrak{d}_p^D$ is reflectable for all $\mathfrak{d}_p^D\in\mathfrak{d}^D$.
\end{definition}

\begin{example}
	SDS is RP because every rotation class of $\mathfrak{d}^{\mathcal{D}_7}$ is reflectable: $\mathfrak{d}_2^{\mathcal{D}_7}\sim_r\mathfrak{d}_7^{\mathcal{D}_7}$, $\mathfrak{d}_3^{\mathcal{D}_7}\sim_r\mathfrak{d}_6^{\mathcal{D}_7}$, and $\mathfrak{d}_4^{\mathcal{D}_7}\sim_r\mathfrak{d}_5^{\mathcal{D}_7}$.
\end{example}

\section{Chromatic Theory}
\label{SEC.Chromatic Theory}

\subsection{Chromatic Sets}

Generalizing diatonicism allows us to generalize chromaticism.

\begin{definition}
	The symbols $\flat$ (\textit{flat}), $\natural$ (\textit{natural}), and $\sharp$ (\textit{sharp}) are called \textit{modifiers}; the set $\mathbb{M}=\{\dots$, $\flat\flat\flat$, $\flat\flat$, $\flat$, $\natural$, $\sharp$, $\sharp\sharp$, $\sharp\sharp\sharp$, $\dots\}$ is called the \textit{modifier set} \cite{Hook}. The function $\mu:\mathbb{Z}\to\mathbb{M}$ defined by the equation $$\mu(x)=\begin{cases} \underbrace{\flat\cdots\flat}_x & \textrm{if } x<0 \\ \natural & \textrm{if } x=0 \\ \underbrace{\sharp\cdots\sharp}_x & \textrm{if } x>0 \end{cases}$$ is called the \textit{modifier map}. We define \textit{modifier addition} by the equation $\mu(x)+\mu(y)=\mu(x+y)$ for all $x,y\in\mathbb{Z}$, and we note that $\mu$ is an isomorphism between $\mathbb{Z}$ and $\mathbb{M}$.
\end{definition}

\begin{definition}
	The \textit{chromatic set} corresponding to a diatonic set $D$ is the set $D^\ast=D\times\mathbb{M}$ of elements called \textit{chromatics}.
\end{definition}

\begin{definition}
	The chromatic set corresponding to $\mathcal{D}_7$  is the \textit{standard chromatic set}, denoted $\mathcal{D}_7^\ast$.
\end{definition}

Given a diatonic $d_i\in D$, an integer $m$, and a modifier $\mu(m)\in\mathbb{M}$, we conventionally denote the chromatic $(d_i,\mu(m))\in D^\ast$ by $(d_i)^m$. When $m$ is known, we can notate chromatics in the form $d_i\mu(m)$, where multiplication of the elements is not implied. For example, we may notate the chromatic $(d_i,\mu(-3))$ by either $(d_i)^{-3}$ or $d_i\flat\flat\flat$, the chromatic $(d_i,\mu(1))$ by either $(d_i)^1$ or $d_i\sharp$, and the chromatic $(d_i,\mu(0))$ by either $(d_i)^0$ or $d_i\natural$. 

\begin{definition}
Chromatics of the form $(d_i)^0$ or $d_i\natural$ are called \textit{natural chromatics}.
\end{definition}

We are careful not to conflate natural chromatics and diatonics. For example, $\textrm{C}\natural\ne\textrm{C}$, since C$\natural$ is a chromatic in the standard chromatic set, whereas C is a diatonic in the standard diatonic set. The distinction between natural chromatics and diatonics is more apparent when using the abstract notation: $(d_i,\mu(0))\in D^\ast$ is clearly a different element than $d_i\in D$.

\subsection{Chromatic Transpositions and Prime Chromatic Orderings}
\label{SSEC.Chromatic Transpositions and Prime Chromatic Orderings}

We formalize the notion of transposing between the chromatics of a chromatic set.

\begin{definition}
	Let $D^\ast$ be a chromatic set and $x_1,x_2\in\mathbb{Z}$; a function $T_{x_1,x_2}:D^\ast\to D^\ast$ defined by the equation $T_{x_1,x_2}((d_i)^m)=(d_{i+x_1})^{m+x_2}$ is called a \textit{chromatic transposition}. If the diatonic cardinality of $D$ is $n$, then the index $i+x_1$ of the diatonic component $d_{i+x_1}\in D$ is understood to be $i+x_1\pmod{n}$.
\end{definition}

\begin{example}
	Let $D_6^\ast=\{d_0,\dots,d_5\}\times\mathbb{M}$ be a chromatic set and suppose we have the chromatic transposition $T_{3,-4}:D_6^\ast\to D_6^\ast$. Then $T_{3,-4}((d_1)^4)=(d_{1+3})^{4+(-4)}=(d_4)^0$ and $T_{3,-4}((d_4)^2)=(d_{4+3})^{2+(-4)}=(d_1)^{-2}$.
\end{example}

We formalize the notion of a modified interval.

\begin{definition}
	Let $D^\ast$ be a chromatic set and $(d_i)^x,(d_j)^y\in D^\ast$ such that $[d_i,d_j]_D$ is a $p$th. A \textit{chromatic interval} is an interval from $(d_i)^x$ to $(d_j)^y$ in $D^\ast$, denoted $[(d_i)^x,(d_j)^y]_{D^\ast}$. We call a chromatic interval $[(d_i)^x,(d_j)^y]_{D^\ast}$ \textit{modified} if and only if $x\ne y$; such a modified interval is called a \textit{$\mu(y-x)$ $p$th}, and the modifier $\mu(y-x)$ is called an \textit{interval modifier}.
\end{definition}

\begin{example}
	Suppose we have the chromatic set $D_9^\ast=\{d_0,\dots,d_8\}\times\mathbb{M}$. Given that the interval from $d_2$ to $d_4$ is a third, the modified interval from $(d_2)^{-1}$ to $(d_4)^1$ is a $\mu(1-(-1))$ third, or simply, a $\sharp\sharp$ third. Given that the interval from $d_5$ to $d_6$ is an eleventh, the modified interval from $d_5\natural$ to $d_6\flat$ is a $\mu(-1-0)$ eleventh, or simply, a $\flat$ eleventh.
\end{example}

We turn to the concept of generating a chromatic ordering from a diatonic space.

\begin{definition}
	Let $D^\ast$ be a chromatic set. Given a diatonic space $\mathfrak{d}_{p,\alpha}^D=\{d_0',\dots,d_{n-1}'\}$, we call the ordering $\dots,d_{n-1}'\flat\flat,d_0'\flat,\dots,d_{n-1}'\flat,d_0'\natural,\dots,d_{n-1}'\natural,d_0'\sharp,\dots,d_{n-1}'\sharp,d_0'\sharp\sharp,\dots$ of $D^\ast$ the \textit{prime chromatic ordering} (\textit{PCO}) generated by $\mathfrak{d}_{p,\alpha}^D$, and when $p$ is known, we call the PCO the \textit{line of $p$ths} generated by $\mathfrak{d}_{p,\alpha}^D$.
\end{definition}

One conventionally assumed PCO of $\mathcal{D}_7^\ast$ is already familiar: $\dots$, B$\flat\flat$, F$\flat$, C$\flat$, G$\flat$, D$\flat$, A$\flat$, E$\flat$, B$\flat$, F$\natural$, C$\natural$, G$\natural$, D$\natural$, A$\natural$, E$\natural$, B$\natural$, F$\sharp$, C$\sharp$, G$\sharp$, D$\sharp$, A$\sharp$, E$\sharp$, B$\sharp$, F$\sharp\sharp$, $\dots$, called the \say{line of fifths} \cite{Temperley}.

\begin{definition}
	The \textit{standard line of fifths} (\textit{standard PCO}) is the line of fifths generated by $\mathfrak{d}_{5,3}^{\mathcal{D}_7}$.
\end{definition}

\begin{definition}
	An \textit{ordered chromatic set} is a chromatic set $D^\ast$ ordered by the line of $p$ths generated by a diatonic space $\mathfrak{d}_{p,\alpha}^{D}$, denoted $\mathfrak{c}_{p,\alpha}^{D^\ast}$.
\end{definition}

\begin{definition}
	The set $\mathfrak{c}_{5,3}^{\mathcal{D}_7^\ast}$ is the \textit{standard ordered chromatic set}.
\end{definition}

The standard ordered chromatic set is \say{the central row of the \textit{Tonnetz} as constructed by Hugo Riemann} \cite{Hook}, so an arbitrary ordered chromatic set could be likened to the central row of a \textit{Tonnetz}-like lattice diagram that accommodates the use of an arbitrary diatonic space \cite{Riemann}.

\begin{definition}
	The terms \textit{grave} and \textit{acute} refer to the leftward and rightward directions of a PCO, respectively. Given an ordered chromatic set $\mathfrak{c}_{p,\alpha}^{D^\ast}$, $(d_i)^{m_1}$ is \textit{graver} than $(d_j)^{m_2}$ (or, equivalently, $(d_j)^{m_2}$ is \textit{acuter} than $(d_i)^{m_1}$) if and only if $(d_i)^{m_1}$ precedes $(d_j)^{m_2}$ in $\mathfrak{c}_{p,\alpha}^{D^\ast}$, denoted $(d_i)^{m_1}\prec(d_j)^{m_2}$ or $(d_j)^{m_2}\succ(d_i)^{m_1}$.
\end{definition}

\begin{example}
	Suppose $D_5=\{d_0,\dots,d_4\}$ is a diatonic set of diatones $d_0=\textrm{A}$, $d_1=\textrm{B}$, $d_2=\textrm{C}$, $d_3=\textrm{D}$, and $d_4=\textrm{E}$. Since $\mathfrak{d}_{3,1}^{D_5}=\{\textrm{B},\textrm{D},\textrm{A},\textrm{C},\textrm{E}\}$, then ordering $D_5^\ast$ with the line of thirds generated by $\mathfrak{d}_{3,1}^{D_5}$ yields $\mathfrak{c}_{3,1}^{D_5^\ast}=\{\dots,$ E$\flat\flat$, B$\flat$, D$\flat$, A$\flat$, C$\flat$, E$\flat$, B$\natural$, D$\natural$, A$\natural$, C$\natural$, E$\natural$, B$\sharp$, D$\sharp$, A$\sharp$, C$\sharp$, E$\sharp$, B$\sharp\sharp$, $\dots\}$. We have D$\flat$ graver B$\natural$ and E$\sharp$ acuter than A$\natural$ in $\mathfrak{c}_{3,1}^{D_5^\ast}$.
\end{example}

Given a diatonic set $D$ with diatonic cardinality $n$, the line of $p$ths ordering of $D^\ast$ generated by $\mathfrak{d}_{p,\alpha}^D$ is a realization of $D^\ast$ as a Generalized Interval System as defined in \textit{Generalized Musical Intervals and Transformations} \cite{Lewin}, where the interval group is $\mathbb{Z}$ such that $\textrm{int}{((d_i')^{m_1},(d_j')^{m_2})}=j-i+(m_2-m_1)n$ \cite{Hook}. If $(d_i')^{m_1}\prec(d_j')^{m_2}$, then $\textrm{int}{((d_i')^{m_1},(d_j')^{m_2})}>0$, and if $(d_i')^{m_1}\succ(d_j')^{m_2}$, then $\textrm{int}{((d_i')^{m_1},(d_j')^{m_2})}<0$. Informally, we can think of the interval as indicating how many \say{steps} we must go in either the leftward or rightward direction to traverse the interval.

\begin{example}
	Using the standard ordered chromatic set, $\textrm{int}{(\textrm{G}\natural,\textrm{B}\flat)}=\textrm{int}{((d_2')^0,(d_6')^{-1})}=6-2+7(-1-0)=-3$ and $\textrm{int}{(\textrm{B}\flat,\textrm{G}\natural)}=3$; that is, to get from G$\natural$ to B$\flat$ in the standard ordered chromatic set, we take $3$ \say{steps} to the left, and to get from B$\flat$ to G$\natural$, we take $3$ \say{steps} to the right.
\end{example}

In sum, Section \ref{SEC.Chromatic Theory} presents chromaticism as the interaction between diatonics and modifiers, where the structure of diatonic spaces informs the structure of ordered chromatic sets. The proposed definitions suggest a non-acoustic motivation for chromaticism: our chromatic theory is built upon modifying diatonic spaces to generate a higher degree of structural complexity within syntactic music theory---a complexity we harness with enharmonic equivalence.

\section{Enharmonic Theory}
\label{SEC.Enharmonic Theory}

\subsection{Preliminaries}
\label{SSEC.Preliminaries}

Relating chromatics and enharmonics is relatively intuitive, particularly for musicians and those who read sheet music. We will show that an enharmonic is an equivalence class containing chromatics that we call \say{enharmonically equivalent.} To obtain such equivalence classes, we must define a special type of equivalence relation such that it partitions a chromatic set into our desired enharmonics.

\begin{figure}[htp]
	\centering
	\includegraphics[scale=0.3]{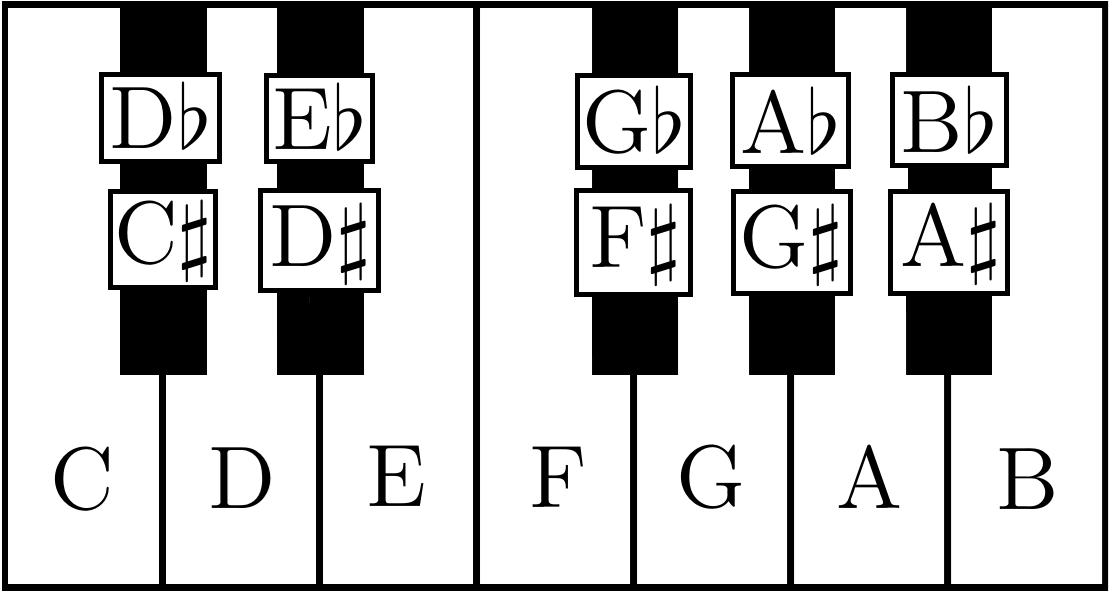}
	\caption{The $12$ standard piano keys (enharmonics).}
	\label{FIG.PianoKeys}
\end{figure}

Defining such an equivalence relation proves to be the main challenge of this section, so before presenting definitions, we first step through the familiarly assumed notion of enharmonic equivalence so that, hopefully, the abstract definitions are more readily understood. It will be helpful to refer to Figure \ref{FIG.PianoKeys} throughout this exploration. Furthermore, we use the terms \say{piano keys} and \say{enharmonics} interchangeably, as a piano key is simply a visual representation of an enharmonic. We will return to this terminology and formalize it in Section \ref{SSEC.Enharmonics and Enharmonic Spaces}, but for the purposes of this exploration, we use the terms informally.

There are $12$ standard piano keys ($7$ white and $5$ black), depicted in Figure \ref{FIG.PianoKeys}. The $7$ white keys appear to correspond to the $7$ diatonics in $\mathcal{D}_7$, and the $5$ black keys are placed between the white keys in a particular arrangement. We show that the arrangement of the black keys with respect to the white keys underlies how we define enharmonic equivalence.

If we let $T_{0,e_i}((d_i)^0)\approx_{\mathcal{E}}(d_{i+1})^0$ denote that $(d_i)^{e_i}$ and $(d_{i+1})^0$ are enharmonically equivalent according to conventional enharmonic equivalence (our use of $\approx_\mathcal{E}$ in Section \ref{SSEC.Preliminaries} is informal, and we formalize enharmonic equivalence in Section \ref{SSEC.Enharmonics and Enharmonic Spaces}), then \begin{align*} T_{0,e_0}(\textrm{C}\natural) &\approx_{\mathcal{E}} \textrm{D}\natural, & T_{0,e_2}(\textrm{E}\natural) &\approx_{\mathcal{E}} \textrm{F}\natural, & T_{0,e_4}(\textrm{G}\natural) &\approx_{\mathcal{E}} \textrm{A}\natural, & T_{0,e_6}(\textrm{B}\natural) &\approx_{\mathcal{E}} \textrm{C}\natural, \\ T_{0,e_1}(\textrm{D}\natural) &\approx_{\mathcal{E}} \textrm{E}\natural, & T_{0,e_3}(\textrm{F}\natural) &\approx_{\mathcal{E}} \textrm{G}\natural, & T_{0,e_5}(\textrm{A}\natural) &\approx_{\mathcal{E}} \textrm{B}\natural, \end{align*} where $e_0=e_1=e_3=e_4=e_5=2,\,e_2=e_6=1$, and $\mathcal{E}=(e_0,\dots,e_6)=(2,2,1,2,2,2,1)$. In other words, we can think of $\mathcal{E}$ as containing the specific instructions to properly assign conventional enharmonic equivalence.

Each element of $\mathcal{E}$ tells us how many black keys to insert between each of the white keys; that is, we insert $e_i-1$ black keys in between the white key containing $(d_i)^0$ and the white key containing $(d_{i+1})^0$. Because we do not want negative or fractional numbers of black keys in between two white keys, we restrict the values of $e_i$ to be strictly positive integers. Although not a relevant constraint within conventional enharmonic equivalence, we will need this property to generalize enharmonic equivalence. Furthermore, we observe that the sum $e_0+\dots+e_6=2+2+1+2+2+2+1=12$ yields the total number of piano keys, white and black combined.

We observe that the equivalences listed previously can all be expressed with respect to C$\natural$, using the values of $\mathcal{E}$: \begin{align*} T_{0,2}(\textrm{C}\natural) &\approx_{\mathcal{E}} \textrm{D}\natural, & T_{0,2+2+1+2}(\textrm{C}\natural) &\approx_{\mathcal{E}} \textrm{G}\natural, & T_{0,2+2+1+2+2+2}(\textrm{C}\natural) &\approx_{\mathcal{E}} \textrm{B}\natural, \\ T_{0,2+2}(\textrm{C}\natural) &\approx_{\mathcal{E}} \textrm{E}\natural, & T_{0,2+2+1+2+2}(\textrm{C}\natural) &\approx_{\mathcal{E}} \textrm{A}\natural, & T_{0,2+2+1+2+2+2+1}(\textrm{C}\natural) &\approx_{\mathcal{E}} \textrm{C}\natural, \\ T_{0,2+2+1}(\textrm{C}\natural) &\approx_{\mathcal{E}} \textrm{F}\natural, \end{align*} where we can more generally state that $T_{0,e_0+\dots+e_{i-1}}((d_0)^0)\approx_{\mathcal{E}}(d_i)^0$. Clearly, the consecutive sums of the elements of $\mathcal{E}$ are useful, so we define a special function $A_{\mathcal{E}}:\mathbb{Z}_7\to\mathbb{Z}$ defined by the equation $$A_{\mathcal{E}}(x)=\begin{cases} 0 & \textrm{if }x=0 \\ \sum_{i\,=\,0}^{x-1}e_i & \textrm{if }x\ne 0 \end{cases},\qquad e_i\in\mathcal{E}\textrm{ for each }i=0,\dots,6.$$ Then $T_{0,A_\mathcal{E}(i)}((d_0)^0)\approx_\mathcal{E}(d_i)^0$, which can also be expressed by $(d_0)^{A_\mathcal{E}(i)}\approx_\mathcal{E}(d_i)^0$. Whereas $\mathcal{E}$ stores the \say{instructions} for assigning conventional enharmonic equivalence, $A_\mathcal{E}$ computationally \say{assigns} conventional enharmonic equivalence.

\begin{figure}[htp]
	\centering
	\includegraphics[scale=0.26]{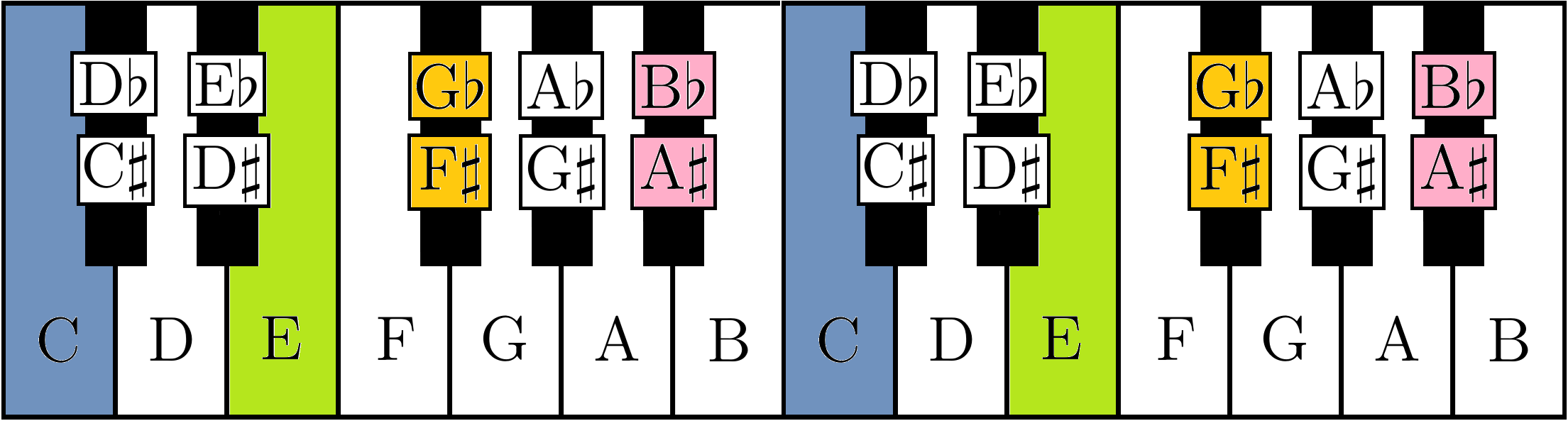}
	\caption{Color-coded examples of enharmonically equivalent piano keys.}
	\label{FIG.PianoKeysTwoOctaves}
\end{figure}

We lastly consider the cyclic nature of a piano keyboard \cite{Fiore,Wright,Tsok}. A piano key must be enharmonically equivalent to itself, and Figure \ref{FIG.PianoKeysTwoOctaves} presents some color-coded examples of this. Since there are $12$ distinct piano keys, then the standard piano keyboard is cyclic modulo $12$, so two chromatics $(d_i)^{m_1}$ and $(d_i)^{m_2}$ are enharmonically equivalent if and only if $m_1\equiv m_2\pmod{12}$.

\subsection{Enharmonic Assignments}
\label{SSEC.Enharmonic Assignments}

Having dissected conventional enharmonic equivalence, we can begin formalizing general enharmonicism.

\begin{definition}
	Let $E=(e_0,\dots,e_{n-1})$ be an $n$-tuple of positive integers, where the indices $i$ of each $e_i$ are the elements of $\mathbb{Z}_n$. The function $A_E:\mathbb{Z}_{n}\to\mathbb{Z}$ defined by the equation $$A_E(x)=\begin{cases} 0 & \textrm{if }x=0 \\ \sum_{i\,=\,0}^{x-1}e_i & \textrm{if }x\ne 0 \end{cases},\qquad e_i\in E\textrm{ for each }i=0,\dots,n-1,$$ is called the \textit{enharmonic assignment} of $E$; the \textit{weight} of an enharmonic assignment $A_E$ is the sum of the elements of $E$.
\end{definition}

\begin{definition}
	The \textit{standard enharmonic assignment} is the enharmonic assignment of $\mathcal{E}=(2,2,1,2,2,2,1)$, denoted $A_\mathcal{E}$.
\end{definition}

We are, of course, not limited to the standard enharmonic assignment.

\begin{example}
	We define a $10$-tuple $E_{10}=(4,2,7,1,1,4,2,6,5,9)$, where the enharmonic assignment of $E_{10}$ is the function $A_{E_{10}}:\mathbb{Z}_{10}\to\mathbb{Z}$ defined by the equation $$A_{E_{10}}(x)=\begin{cases} 0 & \textrm{if }x=0 \\ \sum_{i\,=\,0}^{x-1}e_i & \textrm{if }x\ne 0 \end{cases},\qquad e_i\in E_{10}\textrm{ for each }i=0,\dots,9.$$ The weight of $A_{E_{10}}$ is $4+2+7+1+1+4+2+6+5+9=41$, so we would require $41$ keys to construct a corresponding piano keyboard. Of these $41$ keys, $10$ would be white keys (recall the number of white keys corresponds to the size of the tuple), and the remaining $31$ would be black keys. The black keys would be arranged according to $E_{10}$; there would be $4-1=3$ between the first two white keys, $2-1=1$ between the next two white keys, and so forth. Of course, a physical construction of such a keyboard would be a visual abomination, but what matters is that our definition is general and not restricted to standard assumptions.
\end{example}

\subsection{Enharmonic Sets}
\label{SSEC.Enharmonics and Enharmonic Spaces}

Formalizing enharmonic assignments enables us to formalize enharmonics.

\begin{definition}
	Let $D=\{d_0,\dots,d_{n-1}\}$ be a diatonic set and $D^\ast$ be its corresponding chromatic set; let $E$ be an $n$-tuple of positive integers and $A_E$ be the enharmonic assignment of $E$ of weight $N$. Let $\approx_E$ be the relation defined on $D^\ast$ by $(d_i)^x\approx_E(d_j)^y$ if and only if $x+A_E(i)\equiv y+A_E(j)\pmod{N}$. The relation $\approx_E$ is called the \textit{enharmonic equivalence relation} of $E$, and if $(d_i)^x\approx_E(d_j)^y$, we call $(d_i)^x$ and $(d_j)^y$ \textit{enharmonically equivalent} by $E$.
\end{definition}

\begin{remark}
	That $\approx_E$ is an equivalence relation follows immediately from equivalence modulo $N$.
\end{remark}

\begin{definition}
\label{DEF.Trivial Enharmonic Equivalence Relation}
	Let $D=\{d_0,\dots,d_{n-1}\}$ be a diatonic set and $E=(e_0,\dots,e_{n-1})$ be an $n$-tuple such that $e_i=1$ for each $i=0,\dots,n-1$. Then the enharmonic equivalence relation of $E$ is called the \textit{trivial enharmonic equivalence relation}.
\end{definition}

\begin{definition}
\label{DEF.Enharmonic Set}
	A diatonic set $\mathcal{P}_E=\{\pi_0,\dots,\pi_{N-1}\}$ is called \textit{enharmonic} if and only if there exists a diatonic set $D$ such that the diatones $\pi_i\in\mathcal{P}_E$ are the equivalence classes generated by an enharmonic equivalence relation $\approx_E$ on $D^\ast$. The elements of an enharmonic set $\mathcal{P}_E$ are called \textit{enharmonics}, and the diatonic cardinality of $\mathcal{P}_E$ is called the \textit{enharmonic cardinality}.
\end{definition}

\begin{lemma}
\label{LEM.Diatonic Sets are Trivially Enharmonic}
	Every diatonic set is enharmonic by trivial enharmonic equivalence.
\end{lemma}

Lemma \ref{LEM.Diatonic Sets are Trivially Enharmonic} follows immediately from Definitions \ref{DEF.Trivial Enharmonic Equivalence Relation} and \ref{DEF.Enharmonic Set} and demonstrates why we have greater interest in nontrivial enharmonic equivalence relations, since without nontrivial enharmonic equivalence relations, enharmonicism would be redundant. We may now formalize some common terminology with respect to enharmonics.

\begin{definition}
	A (piano) \textit{key} refers to an enharmonic, a \textit{keyboard} refers to an enharmonic set, a \textit{white key} refers to an enharmonic containing a natural chromatic, and a \textit{black key} refers to an enharmonic that does not contain a natural chromatic.
\end{definition}

\begin{figure}[htp]
	\centering
	\includegraphics[scale=0.34]{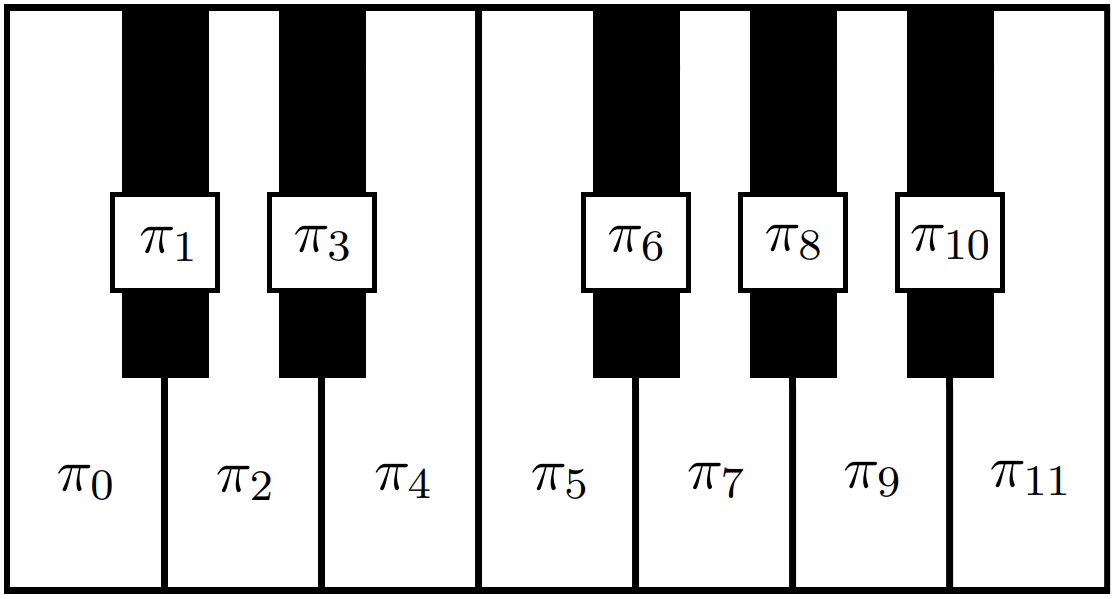}
	\caption{The standard enharmonic set with enharmonics expressed in $\pi$-notation.}
	\label{FIG.Standard Enharmonic Set Pi Notation}
\end{figure}

\begin{definition}
	The enharmonic equivalence relation of $\mathcal{E}=(2,2,1,2,2,2,1)$ on $\mathcal{D}_7^\ast$ is the \textit{standard enharmonic equivalence relation}, denoted $\approx_\mathcal{E}$, the resulting enharmonics are the \textit{standard enharmonics}, which are the elements of the \textit{standard enharmonic set}, denoted $\mathfrak{p}=\{\textrm{C}$, C$\sharp$/D$\flat$, D, D$\sharp$/E$\flat$, E, F, F$\sharp$/G$\flat$, G, G$\sharp$/A$\flat$, A, A$\sharp$/B$\flat$, B$\}$, where the symbol \say{/} denotes interchangeable labels for the same enharmonic (see Figure \ref{FIG.Standard Enharmonic Set Pi Notation}).
\end{definition}

\begin{example}
	We call $\textrm{B}\flat\flat$ enharmonically equivalent to $\textrm{G}\sharp\sharp$; since $\textrm{B}\flat\flat=(d_6)^{-2}$ and $\textrm{G}\sharp\sharp=(d_4)^2$, then $-2+A_\mathcal{E}(6)=9=2+A_\mathcal{E}(4)$.
\end{example}

We clarify some potential acoustic misunderstandings for the standard enharmonics also relevant for arbitrary enharmonic sets. The definitions proposed in this paper do not presuppose equal temperament. Recall from Section \ref{SSEC.Acoustics and Syntactics} that this paper is concerned with solely syntactic considerations. We do, however, address a relevant issue: enharmonically equivalent accidentals often refer to distinct acoustic pitches. With the exception of equal-tempered instruments, musicians produce distinct pitches for C$\sharp$ and D$\flat$, for example. Although a fascinating study within tuning theory, we do not attempt to define any relationship between enharmonic labels and their corresponding pitches in cents to fit within a certain tuning, nor do we suggest any nuance in the use of interchangeable enharmonic labels within this paper. More generally, though the indices of enharmonics are in $\mathbb{Z}_N$, we do not presuppose equal temperament, so we cannot assume transposability; that is, we may not assume equal acoustic intervals when we have keyboards with equal numbers of white and black keys (refer to Section \ref{SSEC.Historical Considerations}).

Concerning notation, we clarify the distinction between the statements $\textrm{E}\sharp\in\textrm{F}$ and $\textrm{E}\sharp\approx_\mathcal{E}\textrm{F}\natural$; we could equivalently express these statements using the notation $(d_2)^1\in\pi_5$ and $(d_2)^1\approx_\mathcal{E}(d_3)^0$, respectively. The first statement tells us $\textrm{E}\sharp=(d_2)^1$ is a chromatic in the enharmonic $\textrm{F}=\pi_5$, and the second statement tells us the chromatics E$\sharp$ and F$\natural$ in $\mathcal{D}_7^\ast$ are enharmonically equivalent by $\mathcal{E}$. We do not conflate F$\natural$ and F, since $(d_3)^0$ is a chromatic, and $\pi_5$ is an enharmonic ($(d_3)^0\in\pi_5$).

We recognize the notational ambiguity between diatonics, chromatics, and enharmonics. For example, the symbol F could refer to either the diatone $d_3$ or the enharmonic $\pi_5$; similarly, the symbol B$\flat$ could refer to either the chromatic $(d_6)^{-1}$ or the enharmonic $\pi_{10}$. Within this paper, context determines meaning, and in situations where context is insufficient, we revert to the unambiguous abstract notation.

\subsection{Major Scales and Enharmonic Intervals}
\label{SSEC.Major Scales and Enharmonic Intervals}

Since enharmonic sets are diatonic sets, we only require a clarification of terminology for \say{enharmonic transposition.}

\begin{definition}
	Given an enharmonic set $\mathcal{P}_E$, a diatonic transposition $t_x:\mathcal{P}_E\to\mathcal{P}_E$ is called an \textit{enharmonic transposition} by $x$, denoted $\tau_x$.
\end{definition}

Except for trivial enharmonic equivalence, it is insufficient to simply call $\tau_1$ the \say{enharmonic step transposition}. Consider standard enharmonicism, where we distinguish between \say{whole steps,} \say{half steps,} and \say{chromatic steps,} such that chromatic steps coincide with half steps. We generalize the notion of these three types of step transpositions to accommodate arbitrary enharmonic equivalence.

\begin{definition}
\label{DEF.Enharmonic Step Types}
	Given an enharmonic set $\mathcal{P}_E$, let $w$ be the largest element of $E$. We call the enharmonic transposition by $w$ a \textit{whole step}, and if $h>1$ is an integer that divides $w$, we call the enharmonic transposition by $w/h$ a \textit{$1/h$ step}. An enharmonic transposition by $1$ is called a \textit{chromatic step}.
\end{definition}

\begin{example}
	Consider $\mathfrak{p}$. The largest element of $\mathcal{E}$ is $w=2$, so $\tau_2$ is the whole step in $\mathfrak{p}$. The only integer $h>1$ that divides $2$ is $h=2$, so $\tau_{2/2}=\tau_1$ is the $1/2$ step (half step) in $\mathfrak{p}$. And of course, the familiar chromatic step $\tau_1$ coincides with the half step in $\mathfrak{p}$. We have that G$\sharp$ is $2$ whole steps above E because $\tau_2(\tau_2(\textrm{E}))=\textrm{G}\sharp$, and B$\flat$ is $3$ half steps below D$\flat$ because $\tau_1(\tau_1(\tau_1(\textrm{B}\flat)))=\textrm{D}\flat$.
\end{example}

\begin{definition}
	A \textit{$\pi_k$ major scale} in an enharmonic set $\mathcal{P}_E$ is an ordered subset $\{\tau_{A_E(i)}(\pi_k)\,|\,i\in\mathbb{Z}_n\}\subseteq\mathcal{P}_E$ such that $\tau_{A_E(i)}(\pi_k)<\tau_{A_E(j)}(\pi_k)$ if and only if $i<j$.
\end{definition}

\begin{example}
	The A major scale in $\mathfrak{p}$ is $\{\tau_{A_\mathcal{E}(i)}(\textrm{A})\,|\,i\in\mathbb{Z}_7\}=\{\textrm{A}$, B, C$\sharp$, D, E, F$\sharp$, G$\sharp\}$, where $\textrm{A}<\textrm{B}$, $\textrm{B}<\textrm{C}\sharp$, and so forth, where $\textrm{G}\sharp\not<\textrm{A}$ in the major scale ordering.
\end{example}

Note that major scales may be thought of as sequences of various steps as defined in Definition \ref{DEF.Enharmonic Step Types}, where in the standard case, any major scale follows the pattern of whole, whole, half, whole, whole, whole. The final familiar half step in the major scale step sequence is \say{missing} because we definitionally start with the first enharmonic of the scale, and we address a simple method for generating the entire sequence in the following remarks. We can take the elements of an $n$-tuple $E=(e_0,\dots,e_{n-1})$ and use them (preserving their order in $E$) to generate a new $n$-tuple $(\tau_{e_0},\dots,\tau_{e_{n-1}})$ that returns the familiar sequence memorized for the major scale in basic music theory.

\begin{example}
	The $7$-tuple $\mathcal{E}=(2,2,1,2,2,2,1)$ generates $(\tau_2,\tau_2,\tau_1,\tau_2,\tau_2,\tau_2,\tau_1)$, and if we let W denote a whole step and H denote a half step, we obtain the familiar $(\textrm{W},\textrm{W},\textrm{H},\textrm{W},\textrm{W},\textrm{W},\textrm{H})$.
\end{example}

\begin{definition}
\label{DEF.Enharmonic Interval}
	Let $\mathcal{P}_E$ be an enharmonic set, $c$ be a nonnegative integer, and $\pi_\ell\in\mathcal{P}_E$ be $c$ chromatic steps above $\pi_k\in\mathcal{P}_E$. A diatonic interval $[\pi_k,\pi_\ell]_{\mathcal{P}_E}$ in $\mathcal{P}_E$ is called an \textit{enharmonic interval}; the \textit{chromatic length} of $[\pi_k,\pi_\ell]_{\mathcal{P}_E}$ is $c$. Suppose $c$ is in the range of $A_E$. We call the interval $[\pi_k,\pi_\ell]_{\mathcal{P}_E}$ \textit{perfect} if and only if $c\ne N-1$ is coprime with $N$; for all other $c\ne 0$, we call the interval \textit{major}, and for the special case $c=0$, we call the interval $[\pi_k,\pi_k]_{\mathcal{P}_E}$ a \textit{perfect unison} of trivial length $0$.
\end{definition}

When defining the notion of a perfect enharmonic interval, we are only interested in the interval lengths in the range of $A_E$ not trivially coprime with $N$. We consider $N-1$ to be trivially coprime with $N$ because the equation $(N-1)x+Ny=1$ has, for all integers $N>1$, the integer solution $x=-1,y=1$, so $\gcd{(N-1,N)}=1$ \cite{Bezout}.

\begin{definition}
\label{DEF.Perfect and Major pth}
	Let $[d_i,d_j]_D$ be a $p$th and $\mathcal{P}_E$ be an enharmonic set. Given $(d_i)^x\in\pi_k$ and $(d_j)^y\in\pi_\ell$ such that $[\pi_k,\pi_\ell]_{\mathcal{P}_E}$ is perfect (or major), we call $[\pi_k,\pi_\ell]_{\mathcal{P}_E}$ a \textit{perfect} (or \textit{major}) \textit{$p$th}.
\end{definition}

\begin{example}
	We call $[\textrm{B}\flat,\textrm{F}]_{\mathfrak{p}}$ a perfect fifth because $\tau_7(\textrm{B}\flat)=\textrm{F}$ and $[\textrm{B},\textrm{F}]_{\mathcal{D}_7}$ is a fifth. It is necessary that $7$ satisfies the following definitional properties: (i) $7$ is in the range of $A_\mathcal{E}$ and (ii) $7\ne 11$ and $12$ are coprime. Similarly, we call $[\textrm{F}\sharp,\textrm{A}\sharp]_{\mathfrak{p}}$ a major third because $\tau_4(\textrm{F}\sharp)=\textrm{A}\sharp$ and $[\textrm{F},\textrm{A}]_{\mathcal{D}_7}$ is a third. Again, (i) $4$ is in the range of $A_\mathcal{E}$ and (ii) $\gcd{(4,12)}\ne 1$.
\end{example}

It will be useful to quickly identify which intervals in a diatonic set are perfect in a corresponding enharmonic set, so we give the following lemma which follows immediately from Definitions \ref{DEF.Enharmonic Interval} and \ref{DEF.Perfect and Major pth}.

\begin{lemma}
\label{LEM.Identifying pths}
	A $p$th in $D$ is perfect in $\mathcal{P}_E$ if and only if $A_E(p-1)\ne N-1$ and $A_E(p-1)$ is coprime with $N$. If a $p$th in $D$ is not perfect in $\mathcal{P}_E$, then it is major in $\mathcal{P}_E$ except for perfect unison.
\end{lemma}

Lemma \ref{LEM.Identifying pths} allows us to quickly verify that the only perfect intervals in $\mathfrak{p}$ are fourths and fifths without needing specific examples.

\begin{example}
	Fourths in $\mathcal{D}_7$ are perfect in $\mathfrak{p}$ because $A_\mathcal{E}(4-1)=5\ne 11$ (and $5$ is coprime with $12$), and fifths in $\mathcal{D}_7$ are also perfect in $\mathfrak{p}$ because $A_\mathcal{E}(5-1)=7\ne 11$ (and $7$ is coprime with $12$). The major intervals in $\mathfrak{p}$ are seconds, thirds, sixths, and sevenths; sixths in $\mathcal{D}_7$ are major in $\mathfrak{p}$ because $\gcd{(A_\mathcal{E}(6-1),12)}=\gcd{(9,12)}\ne 1$.
\end{example}

Note that intervals such as $[\textrm{D},\textrm{F}]_{\mathfrak{p}}$ and $[\textrm{C},\textrm{G}\flat]_{\mathfrak{p}}$ are neither perfect nor major, since $\tau_3(\textrm{D})=\textrm{F}$ and $\tau_6(\textrm{C})=\textrm{G}\flat$, respectively, and neither $3$ nor $6$ satisfy the property of being in the range of $A_\mathcal{E}$. That there can exist enharmonic intervals that are neither perfect nor major motivates the following definition, which generalizes the familiar notions of \say{minor,} \say{augmented,} and \say{diminished} intervals.

\begin{definition}
\label{DEF.Augmented and Diminished Perfect Intervals}
	Given a perfect $p$th $[\pi_k,\pi_\ell]_{\mathcal{P}_E}$, the interval $[\pi_k,\tau_c(\pi_\ell)]_{\mathcal{P}_E}$ is \textit{augmented} if and only if $c>0$ and \textit{diminished} if and only if $c<0$. We inductively define the terminology \textit{augmented} for $c=1$, \textit{double augmented} for $c=2$, \textit{triple augmented} for $c=3$, and so forth. Similarly, we call $[\pi_k,\tau_c(\pi_\ell)]_{\mathcal{P}_E}$ \textit{diminished} for $c=-1$, \textit{double diminished} for $c=-2$, and so forth.
\end{definition}

\begin{example}
	The previously mentioned interval $[\textrm{C},\textrm{G}\flat]_{\mathfrak{p}}$ can now be classified as a diminished fifth. We first identify $[\textrm{C},\textrm{G}]_{\mathfrak{p}}$ as a perfect fifth, since $\tau_7(\textrm{C})=\textrm{G}$ and $[\textrm{C},\textrm{G}]_{\mathcal{D}_7}$ is a fifth. Then we see that $[\textrm{C},\textrm{G}\flat]_{\mathfrak{p}}$ is of the form $[\textrm{C},\tau_{-1}(\textrm{G})]_{\mathfrak{p}}$, and since $c=-1$, we have a diminished fifth.
\end{example}

\begin{definition}
	Given a major $p$th $[\pi_k,\pi_\ell]_{\mathcal{P}_E}$, the interval $[\pi_k,\tau_c(\pi_\ell)]_{\mathcal{P}_E}$ is \textit{augmented} if and only if $c>0$, \textit{minor} if and only if $c=-1$, and \textit{diminished} if and only if $c<-1$. We use the same terminology as in Definition \ref{DEF.Augmented and Diminished Perfect Intervals} with the adjustment of starting diminished at $c=-2$ rather than $c=-1$.
\end{definition}

\begin{example}
	The previously mentioned interval $[\textrm{D},\textrm{F}]_{\mathfrak{p}}$ can now be classified as a minor third. We first identify $[\textrm{D},\textrm{F}\sharp]_{\mathfrak{p}}$ as a major third, since $\tau_4(\textrm{D})=\textrm{F}\sharp$ and $[\textrm{D},\textrm{F}]_{\mathcal{D}_7}$ is a third. Then we see that $[\textrm{D},\textrm{F}]_{\mathfrak{p}}$ is of the form $[\textrm{D},\tau_{-1}(\textrm{F}\sharp)]_{\mathfrak{p}}$, and since $c=-1$, we have a minor third.
\end{example}

We conclude this section with two relevant remarks. Given either a perfect or major $p$th $[\pi_k,\pi_\ell]_{\mathcal{P}_E}$, the upper bound on diminishing the $p$th is $\ell-k$, since $[\pi_k,\tau_{-(\ell-k)}(\pi_\ell)]_{\mathcal{P}_E}=[\pi_k,\pi_k]_{\mathcal{P}_E}$ is a perfect unison, and we do not allow diminishing to continue beyond perfect unison in this paper. Also consider the non-uniqueness of enharmonic interval classification. In addition to being called a diminished fifth, the interval $[\textrm{C},\textrm{G}\flat]_{\mathfrak{p}}$ can also be called an augmented fourth, since $[\textrm{C},\textrm{F}]_{\mathfrak{p}}$ is a perfect fourth, and $[\textrm{C},\tau_1(\textrm{F})]_{\mathfrak{p}}=[\textrm{C},\textrm{F}\sharp]_{\mathfrak{p}}$ (recall F$\sharp$ and G$\flat$ are interchangeable labels for $\pi_6\in\mathfrak{p}$). In fact, there are even more ways to name this interval in both directions of augmenting and diminishing; we could call $[\textrm{C},\textrm{G}\flat]_{\mathfrak{p}}$ something absurd, like a quadruple augmented second or a triple diminished sixth. Such names, of course, are strongly discouraged in practice, and as an informal convention, we seek to classify our intervals according to the minimal amount of augmenting or diminishing, if required (this convention does not imply uniqueness, as both diminished fifth and augmented fourth would be considered appropriate classifications of $[\textrm{C},\textrm{G}\flat]_{\mathfrak{p}}$).

\subsection{Prime Enharmonic Orderings and Enharmonic Systems}
\label{SSEC.Prime Enharmonic Orderings and Enharmonic Systems}

The construction of an enharmonic system from an arbitrary enharmonic set is analogous to the construction of a diatonic system from an arbitrary diatonic set; explanations of motivations in this section are sparse or omitted to avoid redundancy with Sections \ref{SSEC.Diatonic Transpositions and Prime Diatonic Orderings} and \ref{SSEC.Diatonic Systems}.

\begin{definition}
\label{DEF.PEO}
	Given an enharmonic set $\mathcal{P}_E=\{\pi_0,\dots,\pi_{N-1}\}$ with $N>1$ and an integer $p$ such that a $p$th is perfect in $\mathcal{P}_E$, let $\eta_p:\mathcal{P}_E\to\mathcal{P}_E$ be a function defined by the equation $\eta_p(\pi_x)=\pi_{A_E(p-1)x}$. Elements in the range of $\eta_p$ are denoted $\pi_i'$, where $\pi_i'=\eta_p(\pi_i)$, and are called \textit{prime enharmonics}. An ordering $\pi_0',\pi_1',\dots,\pi_{N-1}'$ of $\mathcal{P}_E$ is called a \textit{prime enharmonic ordering} (\textit{PEO}), and if $p$ is known, we call the resulting PEO an \textit{enharmonic circle of perfect $p$ths}.
\end{definition}

Note that Lemma \ref{LEM.Identifying pths} together with the proof given in Section \ref{SSEC.Diatonic Transpositions and Prime Diatonic Orderings} guarantees $\eta_p$ is bijective for all valid $p$.

\begin{lemma}
\label{LEM.Not All PDOs are PEOs}
	All PEOs are PDOs, but not all PDOs are PEOs.
\end{lemma}

Lemma \ref{LEM.Not All PDOs are PEOs} follows from Definitions \ref{DEF.PDO} and \ref{DEF.PEO}. The set of integers $p$ that can be used to generate a PEO is a subset of the set of integers $p$ that can be used to generate a PDO, since we allow $p-1=n-1$ in Definition \ref{DEF.PDO} but exclude $A_E(p-1)=N-1$ in Definition \ref{DEF.PEO}, due to the requirement that the interval be perfect (Definition \ref{DEF.Enharmonic Interval}).

\begin{definition}
	The enharmonic circle of perfect fifths of $\mathfrak{p}$ is the \textit{standard enharmonic circle of perfect fifths} (\textit{standard PEO}).
\end{definition}

Stated explicitly, the standard PEO is C, G, D, A, E, B, F$\sharp$/G$\flat$, C$\sharp$/D$\flat$, G$\sharp$/A$\flat$, D$\sharp$/E$\flat$, A$\sharp$/B$\flat$, F. Of course, conventional introductions to the standard PEO do not present all interchangeable labels, and the labeling C, G, D, A, E, B, F$\sharp$/G$\flat$, D$\flat$, A$\flat$, E$\flat$, B$\flat$, F is preferred. However, such labeling depends on the notion of \say{key signature}, which we do not explicitly formalize in this paper, hence why we present all possible labels of enharmonics to avoid inadvertently suggesting any meaning by arbitrarily preferring one label to another in an arbitrary PEO. In future work, obtaining a generalization of key signature within this framework of definitions would enable us to have a rigorous procedure for \say{preferring} one enharmonic label to another in a PEO.

\begin{definition}
	An \textit{enharmonic space} is a diatonic space of an enharmonic set $\mathcal{P}_E$, denoted $\mathfrak{e}_{p,\alpha}^{\mathcal{P}_E}$, such that the PDO of $\mathfrak{e}_{p,\alpha}^{\mathcal{P}_E}$ is a PEO.
\end{definition}

\begin{definition}
	The diatonic system of an enharmonic set $\mathcal{P}_E=\{\pi_0,\dots,\pi_{N-1}\}$ is called an \textit{enharmonic system}, denoted $\mathfrak{e}^{\mathcal{P}_E}=\{\mathfrak{e}_{p}^{\mathcal{P}_E}\,|\,A_E(p-1)\ne N-1,\,\gcd{(A_E(p-1),N)}=1\}$.
\end{definition}

\begin{definition}
	The enharmonic system of $\mathfrak{p}$ is the \textit{standard enharmonic system} (\textit{SES}). The rotation classes of $\mathfrak{e}^\mathfrak{p}$ are the \textit{standard rotation classes}.
\end{definition}

\begin{definition}
	Let $\mathfrak{e}^{\mathcal{P}_E}$ be an enharmonic system of an enharmonic set $\mathcal{P}_E$, where $|\mathcal{P}_E|=N$. A rotation class $\mathfrak{e}_p^{\mathcal{P}_E}\in\mathfrak{e}^{\mathcal{P}_E}$ is an \textit{enharmonic reflection} of a rotation class $\mathfrak{e}_q^{\mathcal{P}_E}\in\mathfrak{e}^{\mathcal{P}_E}$ if and only if $A_E(p-1)+A_E(q-1)=N$ (that is, if $A_E(p-1)$ and $A_E(q-1)$ are additive inverses in $\mathbb{Z}_N$).
\end{definition}

\begin{example}
	Since there are only $2$ perfect intervals in $\mathfrak{p}$ (fourths and fifths), there are only $2$ standard rotation classes: $\mathfrak{e}^\mathfrak{p}=\{\mathfrak{e}_4^\mathfrak{p},\mathfrak{e}_5^\mathfrak{p}\}$. Figure \ref{FIG.PerfectFourthsandFifths} shows representative enharmonic spaces for the two standard rotation classes. We see that SES is RP because every rotation class of $\mathfrak{e}^\mathfrak{p}$ is reflectable ($\mathfrak{e}_4^\mathfrak{p}\sim_r\mathfrak{e}_5^\mathfrak{p}$), a fact made visually apparent by Figure \ref{FIG.PerfectFourthsandFifths}.
	
	\begin{figure}[htp]
		\begin{center}
			\includegraphics[scale=0.7]{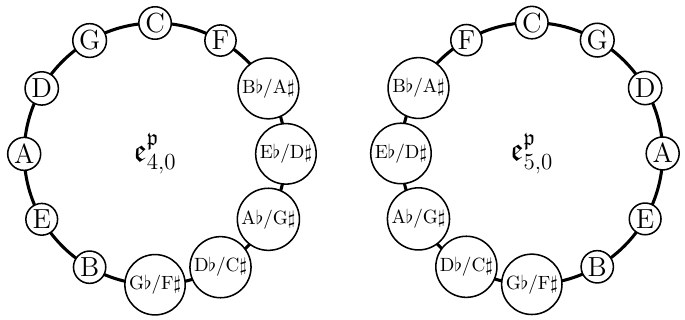}
			\caption{Enharmonic spaces representing their corresponding standard rotation classes.}
			\label{FIG.PerfectFourthsandFifths}
		\end{center}
	\end{figure}
\end{example}

\subsection{Practical Chromatics}
\label{SSEC.Practical Chromatics}

Since the modifier set is isomorphic to the set of integers, $\mathcal{D}_7^\ast$ contains a countably infinite number of chromatics. However, despite there being infinitely many chromatics in each enharmonic of $\mathfrak{p}$, only a small subset of them are used in practice. So, before we consider some examples of nonstandard enharmonic systems, we make some relevant remarks concerning higher modifiers.

It is relatively common for double sharps and double flats to appear in musical scores, where composers use the symbol $\doublesharp$ to denote the modifier $\sharp\sharp$ and $\doubleflat$ to denote $\flat\flat$; going forward, we will denote $\mu(2)=\doublesharp$ and $\mu(-2)=\doubleflat$. However, modifiers beyond $\doublesharp$ and $\doubleflat$ are essentially never used. The theoretical notion of using triple modifiers can be traced back to the Baroque era, where the influential composer Georg Philipp Telemann set forth a framework for understanding triple modifiers \cite{Telemann}. In practice, however, there are extremely limited examples of scores in which the composer deliberately notates a triple sharp or a triple flat.

\begin{figure}[htp]
	\centering
	\includegraphics[scale=0.66]{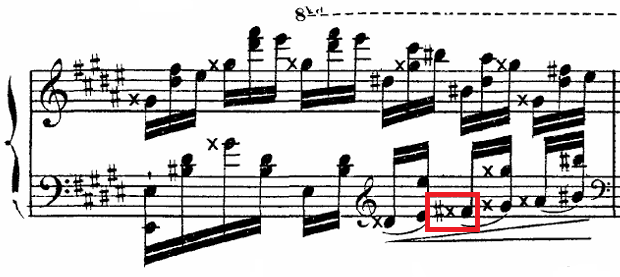}
	\caption{Alkan's triple sharp in measure 291 of his Concerto for Solo Piano, Op. 39, No. 10, Mov. 3.}
	\label{FIG.Alkan's_Triple_Sharp}
\end{figure}

\begin{figure}[htp]
	\centering
	\includegraphics[scale=0.57]{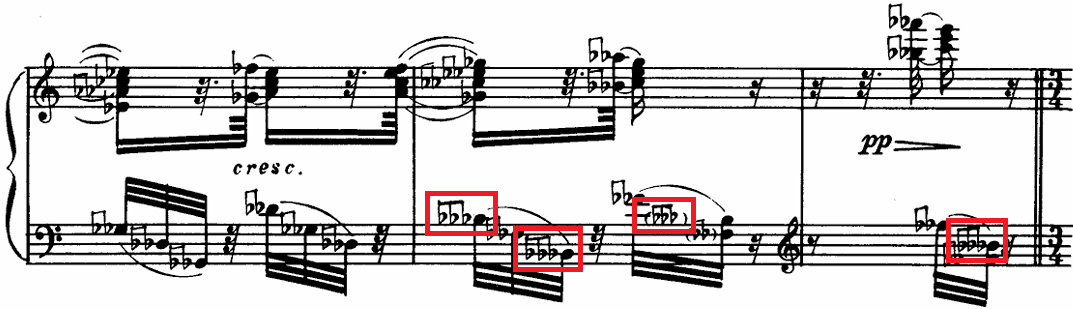}
	\caption{Roslavets' triple flats in measures 151--152 of his Piano Sonata No. 1.}
	\label{FIG.Roslavet's_Triple_Flats}
\end{figure}

A notable example is Charles-Valentin Alkan's Concerto for Solo Piano, Op. 39, No. 10 \cite{Alkan}, which has a triple sharp denoted $\sharp\,\doublesharp$ in measure 291 of the third movement (see Figure \ref{FIG.Alkan's_Triple_Sharp}). Another frequently cited example is Nikolai Andreevich Roslavets' Piano Sonata No. 1 \cite{Roslavets}, which features several triple flats denoted $\flat\flat\flat$ in measures 151--152 (see Figure \ref{FIG.Roslavet's_Triple_Flats}). Although such modifiers do exist, their exceeding rarity causes them to be treated in strictly theoretical contexts for the most part. A major reason is performance and readability---triple modifiers are unnecessarily difficult for musicians to read, especially when enharmonic equivalence can simplify how the score is notated (for example, $\textrm{F}\sharp\,\doublesharp\approx_\mathcal{E}\textrm{G}\sharp$ and $\textrm{B}\flat\flat\flat\approx_\mathcal{E}\textrm{A}\flat$).

Attempting to move beyond triple modifiers is almost never done in standard musical composition and performance, and even in alternate notation systems, there is never practical use of higher modifiers. One of the only composers to tangentially use quadruple modifiers is the modern microtonal pioneer Ben Johnston, whose string quartets sometimes require the use of extended just intonation and an alternate notation system \cite{johnston_notation_2006}. Only within this isolated context does there exist the theoretical possibility of notating quadruple modifiers, but Johnston's music is not practically studied or performed using quadruple modifiers because the alternate notation systems that have been developed to study Johnston's scores are preferred. Thus, although in theory there are a countably infinite number of chromatics in any given enharmonic, in practice, a much smaller subset of those chromatics are actually used to convey music in composition and performance. 

We formalize the notion of what makes a chromatic \say{practical,} not just for $\mathfrak{p}$, but also any arbitrary enharmonic set.

\begin{definition}
\label{DEF.Practical Chromatic}
	Given an enharmonic set $\mathcal{P}_E$ generated by an enharmonic equivalence relation $\approx_E$ on a chromatic set $D^\ast$, let $(d_i)^m\in D^\ast$. The chromatic $(d_i)^m$ is \textit{practical} for $\mathcal{P}_E$ if and only if $-e_{i-1}\le m\le e_i$, where $e_{i-1}$ and $e_i$ are elements of $E$; otherwise, $(d_i)^m$ is \textit{impractical}.
\end{definition}

Since there exist infinitely many nonstandard enharmonic spaces that have practical chromatics with modifier components less than $-2$ or greater than $2$, we will denote such modifiers by $\mu(m)=\sharp^m$ and $\mu(-m)=\flat^m$ for integers $m>2$.

\begin{example}
	To find the practical chromatics with diatonic component E, we recall $\textrm{E}=(d_2)^m\in\mathcal{D}_7^\ast$, so $(d_2)^m$ is practical if and only if $-e_{2-1}\le m\le e_2$, where $e_1=2$ and $e_2=1$ in $\mathcal{E}$; thus, the chromatics E$\doubleflat$, E$\flat$, E$\natural$, and E$\sharp$ are practical, whereas chromatics like E$\flat^5$ and E$\,\doublesharp$ are impractical for $\mathfrak{p}$. We repeat this process for each diatonic component C through B to identify the practical chromatics for $\mathfrak{p}$ (see Table \ref{TAB.Practical Chromatics for p}).
\end{example}

More importantly, Definition \ref{DEF.Practical Chromatic} gives us a rigorous way to restrict the amount of \say{useful} chromatics in general, so that we do not have to rely on arbitrarily selecting which chromatics to include or not in unfamiliar nonstandard enharmonic sets.

\begin{table}[htp]{
	\begin{center}
	\begin{tabular}{c|c}
		Diatonic & Practical Chromatics \\ \hline \\[-1em]
		C & C$\flat$, C$\natural$, C$\sharp$, C$\,\doublesharp$ \\ \\[-1em]
		D & D$\doubleflat$, D$\flat$, D$\natural$, D$\sharp$, D$\,\doublesharp$ \\ \\[-1em]
		E & E$\doubleflat$, E$\flat$, E$\natural$, E$\sharp$ \\ \\[-1em]
		F & F$\flat$, F$\natural$, F$\sharp$, F$\,\doublesharp$ \\ \\[-1em]
		G & G$\doubleflat$, G$\flat$, G$\natural$, G$\sharp$, G$\,\doublesharp$ \\ \\[-1em]
		A & A$\doubleflat$, A$\flat$, A$\natural$, A$\sharp$, A$\,\doublesharp$ \\ \\[-1em]
		B & B$\doubleflat$, B$\flat$, B$\natural$, B$\sharp$ \\
	\end{tabular}
	\end{center}}
	\caption{The practical chromatics for $\mathfrak{p}$.}
	\label{TAB.Practical Chromatics for p}
\end{table}

\subsection{Nonstandard Enharmonic Systems}
\label{SSEC.Nonstandard Enharmonic Systems}

As already discussed in Section \ref{SSEC.Historical Considerations}, standard enharmonicism arises from primarily acoustic considerations, so we hope to show that there are also purely syntactic motivations that lead to preferring SES over arbitrary enharmonic systems. We consider it beautiful that SES can be thought of as the inevitable central construction that unites the acoustic and syntactic branches of music theory. We concede the irony of expending the effort to generalize enharmonicism for the purpose of demonstrating why we prefer the existing SES over the enharmonic systems thus generalized, but in this paper's defense, without the generalized definitions of otherwise familiar concepts, the underlying intuition motivating a syntactic preference of standard enharmonicism is veiled.

Furthermore, considering how musicians continue to poke at the harmonic norms established by standard enharmonicism, we also hope this paper offers a useful syntax for composers to express their musical experimentation in a universally consistent manner, thus making it easier to standardize the transmission of harmonic ideas within nonstandard enharmonic systems. We believe such clarity in communication between musicians would benefit modern composers and the continued expansion of nonstandard harmonic inquiry. Our sentiment here is far from novel; for example, Howe gives a mapping of $\mathfrak{p}_{19}$ (defined below in Example \ref{EX.p19}) onto the traditional five-line staff used for $\mathfrak{p}$, but since $\mathfrak{p}$ and $\mathfrak{p}_{19}$ both originate from the standard diatonic set, such a mapping is logical and quite readable in practice \cite{Howe}.

We show our enharmonic theory is versatile and consistent with existing enharmonic studies. We give examples to show how our proposed definitions can produce nonstandard enharmonic systems that can be analyzed and conveyed using relatively intuitive notation (from the backdrop of standard enharmonicism and syntax). The first example is well known, but the second is arbitrarily constructed for the sake of demonstrating generality and motivating our second important property of useful enharmonic systems, maximal evenness.

\begin{example}
\label{EX.p19}
	We define a nonstandard $7$-tuple $\mathcal{E}'=(3,3,2,3,3,3,2)$ and use the enharmonic equivalence relation $\approx_{\mathcal{E}'}$ to partition $\mathcal{D}_7^\ast$ into $19$ enharmonics labeled $\pi_0=\textrm{C}$, $\pi_1=\textrm{C}\sharp$, $\pi_2=\textrm{D}\flat$, $\pi_3=\textrm{D}$, $\pi_4=\textrm{D}\sharp$, $\pi_5=\textrm{E}\flat$, $\pi_6=\textrm{E}$, $\pi_7=\textrm{E}\sharp/\textrm{F}\flat$, $\pi_8=\textrm{F}$, $\pi_9=\textrm{F}\sharp$, $\pi_{10}=\textrm{G}\flat$, $\pi_{11}=\textrm{G}$, $\pi_{12}=\textrm{G}\sharp$, $\pi_{13}=\textrm{A}\flat$, $\pi_{14}=\textrm{A}$, $\pi_{15}=\textrm{A}\sharp$, $\pi_{16}=\textrm{B}\flat$, $\pi_{17}=\textrm{B}$, and $\pi_{18}=\textrm{B}\sharp/\textrm{C}\flat$. We denote the set of these enharmonics by $\mathfrak{p}_{19}$ (see Figure \ref{FIG.P19KeyboardLabeled}) and give the practical chromatics for $\mathfrak{p}_{19}$ in Table \ref{TAB.Practical Chromatics for p19}.
	
	\begin{figure}[htp]
		\centering
		\includegraphics[scale=0.53]{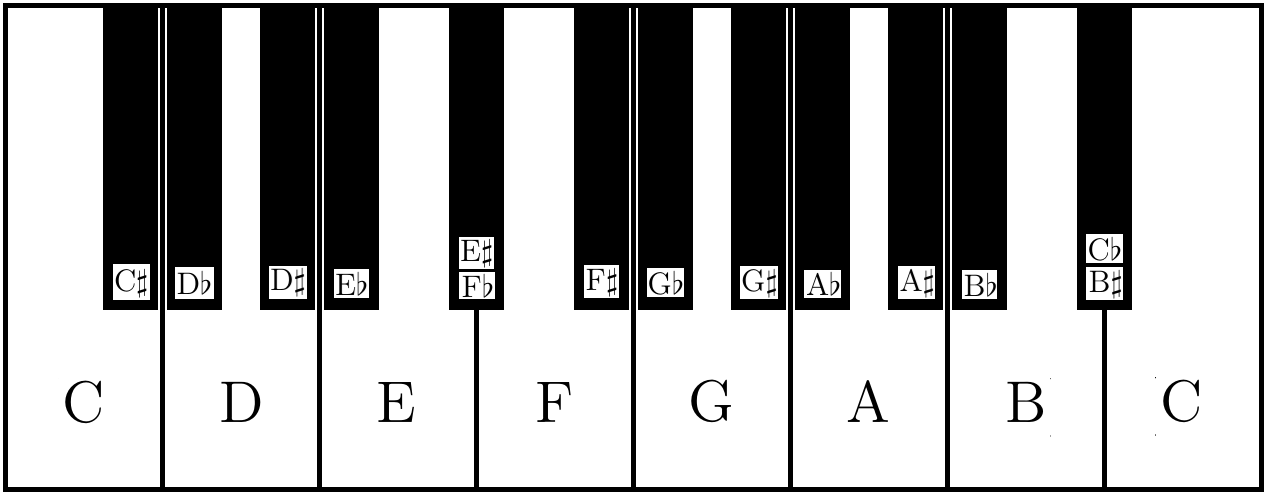}
		\caption{The keyboard representing the nonstandard enharmonic set $\mathfrak{p}_{19}$.}
		\label{FIG.P19KeyboardLabeled}
	\end{figure}
	
	\begin{table}[htp]{
		\begin{center}
		\begin{tabular}{c|c}
				Diatonic & Practical Chromatics \\ \hline \\[-1em]
				C & C$\doubleflat$, C$\flat$, C$\natural$, C$\sharp$, C$\,\doublesharp$, C$\sharp^3$ \\ \\[-1em]
				D & D$\flat^3$, D$\doubleflat$, D$\flat$, D$\natural$, D$\sharp$, D$\,\doublesharp$, D$\sharp^3$ \\ \\[-1em]
				E & E$\flat^3$, E$\doubleflat$, E$\flat$, E$\natural$, E$\sharp$, E$\,\doublesharp$ \\ \\[-1em]
				F & F$\doubleflat$, F$\flat$, F$\natural$, F$\sharp$, F$\,\doublesharp$, F$\sharp^3$ \\ \\[-1em]
				G & G$\flat^3$, G$\doubleflat$, G$\flat$, G$\natural$, G$\sharp$, G$\,\doublesharp$, G$\sharp^3$ \\ \\[-1em]
				A & A$\flat^3$, A$\doubleflat$, A$\flat$, A$\natural$, A$\sharp$, A$\,\doublesharp$, A$\sharp^3$ \\ \\[-1em]
				B & B$\flat^3$, B$\doubleflat$, B$\flat$, B$\natural$, B$\sharp$, B$\,\doublesharp$ \\
		\end{tabular}
		\end{center}}
		\caption{The practical chromatics for $\mathfrak{p}_{19}$.}
		\label{TAB.Practical Chromatics for p19}
	\end{table}
	
	To obtain the enharmonic system of $\mathfrak{p}_{19}$, we must identify the perfect intervals in $\mathfrak{p}_{19}$, but we observe that the enharmonic cardinality is itself prime, meaning that every valid integer $p$ produces a perfect $p$th in $\mathfrak{p}_{19}$; thus, seconds, thirds, fourths, fifths, sixths, and sevenths are perfect in $\mathfrak{p}_{19}$ (note that although sevenths are major in $\mathfrak{p}$, they are perfect in $\mathfrak{p}_{19}$ because $A_{\mathcal{E}'}(7-1)=17\ne 19-1$), so we have the enharmonic system $\mathfrak{e}^{\mathfrak{p}_{19}}=\{\mathfrak{e}_2^{\mathfrak{p}_{19}},\mathfrak{e}_3^{\mathfrak{p}_{19}},\mathfrak{e}_4^{\mathfrak{p}_{19}},\mathfrak{e}_5^{\mathfrak{p}_{19}},\mathfrak{e}_6^{\mathfrak{p}_{19}},\mathfrak{e}_7^{\mathfrak{p}_{19}}\}$.
	
	Computing the chromatic lengths of each perfect interval in $\mathfrak{p}_{19}$ yields the subset $\{3,6,8,11,14,17\}$ of the range of $A_{\mathcal{E}'}$, and of these integers, only $8$ and $11$ form a pair of additive inverses in $\mathbb{Z}_{19}$. Thus, it is straightforward to recognize $\mathfrak{e}^{\mathfrak{p}_{19}}$ is not RP because $\mathfrak{e}_2^{\mathfrak{p}_{19}}$, $\mathfrak{e}_3^{\mathfrak{p}_{19}}$, $\mathfrak{e}_6^{\mathfrak{p}_{19}}$, and $\mathfrak{e}_7^{\mathfrak{p}_{19}}$ are not reflectable in $\mathfrak{p}_{19}$. It is intriguing to observe that, although $\mathfrak{e}^{\mathfrak{p}_{19}}$ is not RP, we still have $\mathfrak{e}_4^{\mathfrak{p}_{19}}\sim_r\mathfrak{e}_5^{\mathfrak{p}_{19}}$, so in some sense, the abstract structure of $\mathfrak{e}^{\mathfrak{p}_{19}}$ behaves \say{nicely} within the circle of perfect fifths and its corresponding opposite circle of perfect fourths (see Figure \ref{FIG.PerfectFourthsandFifthsinp19}).
	
	\begin{figure}[htp]
		\centering
		\includegraphics[scale=0.75]{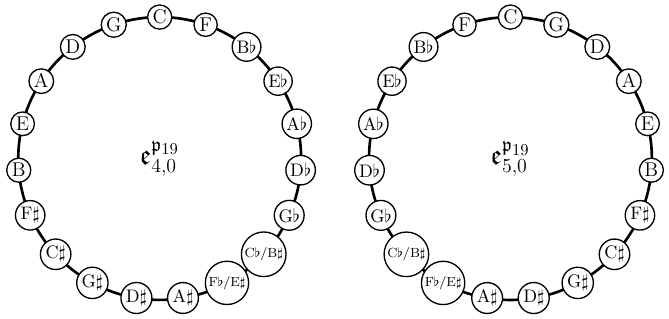}
		\caption{Some representative enharmonic spaces corresponding to their reflectable rotation classes in $\mathfrak{e}^{\mathfrak{p}_{19}}$.}
		\label{FIG.PerfectFourthsandFifthsinp19}
	\end{figure}
	
\end{example}

Note, for example, that although $\textrm{G}\sharp$ and $\textrm{A}\flat$ are the same enharmonic in $\mathfrak{p}$, they are distinct enharmonics in $\mathfrak{p}_{19}$; this difference occurs because $\approx_\mathcal{E}$ and $\approx_{\mathcal{E}'}$ generate different enharmonics on $\mathcal{D}_7^\ast$. Also note that Figure \ref{FIG.P19KeyboardLabeled} does not define an acoustic ordering of $\mathfrak{p}_{19}$; that is, G$\sharp$ in the white key G being visually to the left of A$\flat$ in the white key A does not imply that the pitch of G$\sharp$ is lower than the pitch of A$\flat$. In fact, if we were addressing acoustic considerations in this paper, the pitch used for G$\sharp$ could be higher than the pitch used for A$\flat$ (recall Figure \ref{FIG.Standard Enharmonic Set Pi Notation} and its related remarks about equal-temperament). We clarify that, according to this paper's definitions, neither diatones nor accidentals must preserve pitch order, and consequently, such a remark also applies to the white and black keys of a keyboard. For all keyboard figures (enharmonic sets), we do not presuppose any pitch order and treat such acoustic considerations as outside the scope of this paper.

Among nonstandard enharmonic systems, $\mathfrak{p}_{19}$ is well known and has already been studied in great length. Mandelbaum's 1961 dissertation argues that $\mathfrak{p}_{19}$ is an independently coherent system and not merely a modification of $\mathfrak{p}$ \cite{Mandelbaum}. Mandelbaum's work suggests that $\mathfrak{p}_{19}$ is attractive primarily because it can be generated by the fifth. The notion of a \say{fifth} in Mandelbaum's work arises from the tempered fifth (in $19$-TET) being $11/19$ of the octave (roughly $694.74$ cents), so it is clear his construction of the circle of fifths in $\mathfrak{p}_{19}$ is acoustically motivated (there is a slight discrepancy of about $7$ cents from the Pythagorean fifth, but Mandelbaum argues this difference to be musically acceptable). Not only does Mandelbaum show that $\mathfrak{p}_{19}$ is generated by acoustic fifths, but his observation that the fifth in $\mathfrak{p}_{19}$ consists of $11$ chromatic steps leads to the structural remark that, since $11$ and $19$ are coprime, then $\mathfrak{p}_{19}$ is generated by syntactic fifths as well.

Our definitions do not alter or differ from Mandelbaum's work on $\mathfrak{p}_{19}$, and we hope the enharmonic theory proposed offers the language to abstractly argue for why we tend to prefer certain enharmonic systems over others (in particular, enharmonic systems with a reflectable class of fifths). For example, we could further Mandelbaum's argument for the recognition of $\mathfrak{p}_{19}$ as an independent enharmonic system by pointing to the fact that $\mathfrak{e}_5^{\mathfrak{p}_{19}}$ is reflectable. Even though our enharmonic theory is built independent of pitch considerations, we still arrive to the shared conclusion that $\mathfrak{p}_{19}$ has some musically \say{useful} structure, if we presume that reflectability corresponds to musical practicality.

Gamer's article, \say{Some Combinational Resources of Equal-Tempered Systems} analyzes $\mathfrak{p}_{19}$ with respect to what he calls pitch \say{collections,} \say{transpositions} (as understood in Babbitt's article \say{The Structure and Function of Musical Theory: I} \cite{BabbittArticle}), and \say{interval structures.} \cite{Gamer}. Gamer observed the $\mathbb{Z}_{19}$ group structure of $\mathfrak{p}_{19}$ and knew of the circle of fifths ordering in $\mathfrak{p}_{19}$. His article addresses other ideas we do not address here, as the domain of his study does not intersect with the abstract treatment of enharmonics. His work gives a combinatorial framework that allows us to determine whether a collection or transposition is structurally distinct in $\mathfrak{p}_{19}$, and the article presupposes acoustic ideas like equal temperament to address issues related to microtonal composition, which, at the time of its publication, was increasing in popularity in the West \cite{Gamer}. Gamer's work ultimately provides musicians with a basis for an acoustic $19$-tone theory involving scales, keys, transpositions, and harmony. As with Mandelbaum, we again reach the following remark: though our definitions are divorced from acoustic considerations, our example of $\mathfrak{p}_{19}$ shows how our abstract enharmonic theory accommodates nonstandard structures, even those originating from purely acoustic motivations.

\begin{figure}[htp]
	\centering
	\includegraphics[scale=0.45]{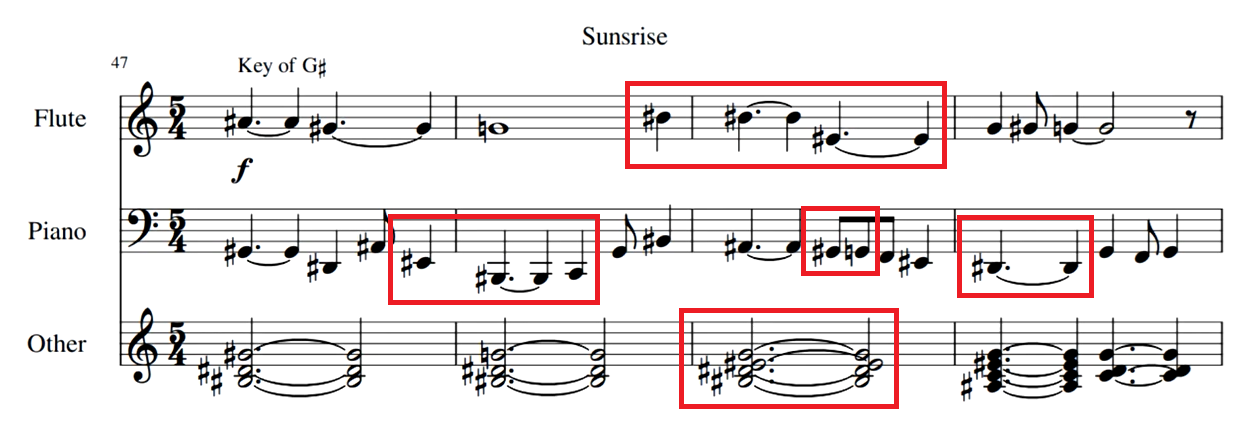}
	\caption{An example of a score composed for $\mathfrak{p}_{19}$.}
	\label{FIG.SunsriseExample}
\end{figure}

Though our purpose in this example is not to analyze scores of music in $\mathfrak{p}_{19}$, we offer some suggestions before continuing: Hook recommends Tracks 1 and 13 of Blackwood's album \textit{Microtonal} \cite{Blackwood} in his article on enharmonicism \cite{Hook}. Neil Haverstick's live performance of \say{Birdwalk} is an example of exploiting $\mathfrak{p}_{19}$ within the blues genre; throughout the composition, enharmonic distinctions between conventionally equivalent enharmonics are melodically employed. Haverstick discusses his notation and composition process for \say{Birdwalk} in his article, \say{Making Microtonal Music Using the $19$-tone Equal Tempered System} \cite{HaverstickBirdwalk}, and the composition receives a more thorough treatment on page 186 of Gann's text, \textit{The Arithmetic of Listening: Tuning Theory and History for the Impractical Musician} \cite{Gann}. Haverstick's piece, \say{Mysteries}, is another of his compositions exploring the sounds of $\mathfrak{p}_{19}$ \cite{HaverstickMysteries}. See Figure \ref{FIG.SunsriseExample} for an example from a score entitled \say{Sunsrise}, also composed using $19$-TET, where the composer (who goes by the artist name Supahstar Saga) makes use of the additional enharmonics $\textrm{B}\sharp\ne\textrm{C}$ and $\textrm{E}\sharp\ne\textrm{F}$ otherwise unavailable in $\mathfrak{p}$ \cite{Saga}. Note in particular how naturally the standard five-line staff system accommodates $\mathfrak{p}_{19}$, which has aided in the transmission of musical compositions in $\mathfrak{p}_{19}$.

\begin{example}
\label{EX.BadEnharmonicSystem}
	Let $D_5=\{d_0,d_1,d_2,d_3,d_4\}$ be a diatonic set of diatones $d_0=\textrm{H}$, $d_1=\textrm{I}$, $d_2=\textrm{J}$, $d_3=\textrm{K}$, and $d_4=\textrm{L}$. Let $E=(3,1,2,4,1)$ and use the enharmonic equivalence relation $\approx_E$ to partition $D_5$ into $11$ enharmonics labeled $\pi_0=\textrm{H}$, $\pi_1=\textrm{H}\sharp$, $\pi_2=\textrm{I}\flat$, $\pi_3=\textrm{I}$, $\pi_4=\textrm{J}$, $\pi_5=\textrm{J}\sharp/\textrm{K}\flat$, $\pi_6=\textrm{K}$, $\pi_7=\textrm{K}\sharp$, $\pi_8=\textrm{K}\,\doublesharp/\textrm{L}\doubleflat$, $\pi_9=\textrm{L}\flat$, and $\pi_{10}=\textrm{L}$. We denote the set of these enharmonics by $\mathcal{P}_E$ (see Figure \ref{FIG.P11KeyboardLabeled}) and give the practical chromatics for $\mathcal{P}_E$ in Table \ref{TAB.Practical Chromatics for P_E}.
	
	\begin{figure}[htp]
		\centering
		\includegraphics[scale=0.55]{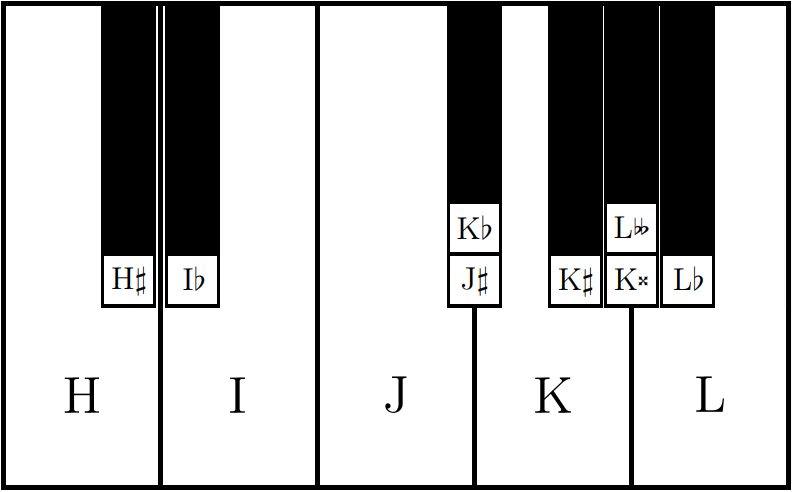}
		\caption{The keyboard representing the nonstandard enharmonic set $\mathcal{P}_{E}$.}
		\label{FIG.P11KeyboardLabeled}
	\end{figure}
	
	\begin{table}[htp]{
		\begin{center}
		\begin{tabular}{c|c}
				Diatonic & Practical Chromatics \\ \hline \\[-1em]
				H & H$\flat$, H$\natural$, H$\sharp$, H$\,\doublesharp$, H$\sharp^3$ \\ \\[-1em]
				I & I$\flat^3$, I$\doubleflat$, I$\flat$, I$\natural$, I$\sharp$ \\ \\[-1em]
				J & J$\flat$, J$\natural$, J$\sharp$, J$\,\doublesharp$ \\ \\[-1em]
				K & K$\doubleflat$, K$\flat$, K$\natural$, K$\sharp$, K$\,\doublesharp$, K$\sharp^3$, K$\sharp^4$ \\ \\[-1em]
				L & L$\flat^4$, L$\flat^3$, L$\doubleflat$, L$\flat$, L$\natural$, L$\sharp$ \\
		\end{tabular}
		\end{center}}
		\caption{The practical chromatics for $\mathcal{P}_{E}$.}
		\label{TAB.Practical Chromatics for P_E}
	\end{table}
	
	We identify the perfect intervals in $\mathcal{P}_E$: seconds, thirds, and fourths (note fifths are major in $\mathcal{P}_E$, rather than perfect, since $A_E(5-1)=10=11-1$). Thus, the enharmonic system of $\mathcal{P}_E$ is $\mathfrak{e}^{\mathcal{P}_E}=\{\mathfrak{e}_2^{\mathcal{P}_E},\mathfrak{e}_3^{\mathcal{P}_E},\mathfrak{e}_4^{\mathcal{P}_E}\}$. Not only does $\mathfrak{e}^{\mathcal{P}_E}$ fail to be RP, but it also fails to contain a single reflectable rotation class: the subset $\{A_E(p-1)\,|\,p=2,3,4\}=\{3,4,6\}$ of the range of $A_E$ contains no pairs of additive inverses in $\mathbb{Z}_{11}$.
\end{example}

From the keyboard representation alone, it is visually apparent that $\mathcal{P}_E$ is irregular and awkwardly defined (caused by the irregularity of $E$). Though not expressly rigorous, we find it reasonable that Example \ref{EX.BadEnharmonicSystem} demonstrates not every enharmonic system is equally conducive for musical application, which motivates the natural question: can we articulate any properties that help us determine the usefulness of an arbitrary enharmonic system? Our formalization of RP (Definition \ref{DEF.Reflection Property}) is one contribution toward answering this question, where we propose RP enharmonic systems are preferred over non-RP enharmonic systems. Example \ref{EX.BadEnharmonicSystem} also motivates the incorporation of Clough and Douthett's notion of \say{maximal evenness} as a second property for determining the usefulness of an enharmonic system \cite{Clough_Douthett}. We discuss our incorporation of Clough and Douthett's notion of \say{ME sets} in Section \ref{SSEC.Maximal Evenness} and ultimately argue that reflectability and maximal evenness are the two basic properties of any musically viable enharmonic system.

\subsection{Maximal Evenness}
\label{SSEC.Maximal Evenness}

Clough and Douthett give the following intuition for their concept of maximal evenness: a maximally even set is \say{a set whose elements are distributed as evenly as possible around the chromatic circle,} where a \say{chromatic circle} in their article is equivalent to our notion of an enharmonic space \cite{Clough_Douthett}. Since the terminology in their paper differs from ours, we clarify some semantic differences. For the following, the terms in quotation marks are used by Clough and Douthett, and the terms in parentheses are their equivalents in our paper: \say{chromatic universe} (\textit{enharmonic set}), \say{chromatic cardinality} (\textit{enharmonic cardinality}), \say{clen} (\textit{chromatic length}), and \say{dlen} (\textit{diatonic length}).

Clough and Douthett define the following concept key to understanding maximal evenness: \say{The spectrum of a dlen is the set of clens corresponding to that particular dlen,} and we write \say{$\langle I\rangle=\{i_1,i_2,\dots\}$ to indicate that the spectrum of dlen $I$ is $\{i_1,i_2,\dots\}$} \cite{Clough_Douthett}. We translate this concept into the context of this paper's definitions in the following.

\begin{definition}
	Let $D=\{d_0,\dots,d_{n-1}\}$ be a diatonic set and $\approx_E$ be an enharmonic equivalence relation on $D^\ast$ such that $E=(e_0,\dots,e_{n-1})$ and $\mathcal{P}_E=\{\pi_0,\dots,\pi_{N-1}\}$ is an enharmonic set. The \textit{spectrum} of a diatonic $d_i\in D$ is the subset $\langle d_i\rangle=\{\tau_x(\pi_{A_E(i)})\,|\,x=0,\dots,e_i-1\}$ of $\mathcal{P}_E$.
\end{definition}

\begin{lemma}
\label{LEM.Spectrums Unique}
	The spectrum of a diatonic is unique and $|\langle d_i\rangle|=e_i$ for all $i\in\mathbb{Z}_n$.
\end{lemma}

Lemma \ref{LEM.Spectrums Unique} follows quickly from former definitions, so we omit a tedious proof, which does not contribute to the purpose of this discussion.

\begin{example}
	The spectra of the standard diatones in $\mathfrak{p}$ are: \begin{align*} \langle\textrm{C}\rangle &= \{\textrm{C},\textrm{C}\sharp/\textrm{D}\flat\}, & \langle\textrm{F}\rangle &= \{\textrm{F},\textrm{F}\sharp/\textrm{G}\flat\}, & \langle\textrm{A}\rangle &= \{\textrm{A},\textrm{A}\sharp/\textrm{B}\flat\}, \\ \langle\textrm{D}\rangle &= \{\textrm{D},\textrm{D}\sharp/\textrm{E}\flat\}, & \langle\textrm{G}\rangle &= \{\textrm{G},\textrm{G}\sharp/\textrm{A}\flat\}, & \langle\textrm{B}\rangle &= \{\textrm{B}\}. \\ \langle\textrm{E}\rangle &= \{\textrm{E}\}, \end{align*}
\end{example}

Clough and Douthett define an enharmonic set to be \say{maximally even} if and only if $1\le|\langle d_i\rangle|\le 2$ for all $d_i$. This definition works well within the confines of standard enharmonicism, but we require a more general definition for our purposes. We interpret Clough and Douthett's definition of maximally even to state that an enharmonic set is maximally even so long as the difference between $|\langle d_i\rangle|$ and $|\langle d_i\rangle|$ is at most $1$ for all $d_i$ and $d_j$, so we formalize our use of the term \say{maximally even} in the following definition.

\begin{definition}
	Let $\mathcal{P}_E$ be an enharmonic set and denote the smallest element of $E$ by $e_{\textrm{min}}$. Then $\mathcal{P}_E$ is \textit{maximally even} if and only if $e_{\textrm{min}}\le|\langle d_i\rangle|\le e_{\textrm{min}}+1$ for all spectra $\langle d_i\rangle$ of $\mathcal{P}_E$.
\end{definition}

\begin{definition}
	An enharmonic system $\mathfrak{e}^{\mathcal{P}_E}$ is said to be \textit{maximally even} (\textit{ME}) if and only if $\mathcal{P}_E$ is maximally even.
\end{definition}

\begin{example}
	SES is ME. The smallest element of $\mathcal{E}$ is $1$, and $1\le |\langle d_i\rangle|\le 2$ for all spectra $\langle d_i\rangle$ of $\mathfrak{p}$.
\end{example}

\begin{example}
	The enharmonic system $\mathfrak{e}^{\mathfrak{p}_{19}}$ is ME. The smallest element of $\mathcal{E}'$ is $2$, and $2\le |\langle d_i\rangle|\le 3$ for all spectra $\langle d_i\rangle$ of $\mathfrak{p}_{19}$.
\end{example}

\begin{example}
	The enharmonic system $\mathfrak{e}^{\mathcal{P}_E}$ defined in Example \ref{EX.BadEnharmonicSystem} is not ME. The smallest element of $E$ is $1$, and there exist spectra such that the inequality $1\le |\langle d_i\rangle|\le 2$ does not hold; for example, $|\langle\textrm{K}\rangle|=|\{\textrm{K},\textrm{K}\sharp,\textrm{K}\,\doublesharp/\textrm{L}\doubleflat,\textrm{L}\flat\}|=4$.
\end{example}

\begin{definition}
	An enharmonic system that is RP and ME is called an \textit{RP-ME enharmonic system}.
\end{definition}

\begin{example}
	Among the three enharmonic systems defined in this paper, only SES is RP-ME. The enharmonic system $\mathfrak{e}^{\mathfrak{p}_{19}}$ is ME but fails to be RP, and the enharmonic system defined in Example \ref{EX.BadEnharmonicSystem} is neither ME nor RP.
\end{example}

\section{Concluding Remarks}
\label{SEC.Concluding Remarks}

We suggest that, within a theory of general enharmonicism, we prefer RP-ME enharmonic systems for musical application and analysis. We hope our proposed enharmonic theory sufficiently demonstrates that it is possible to reach conclusions preferring SES over previously considered nonstandard enharmonic systems from an abstract perspective, fully divorced from pitch considerations. Given the scope of this paper, many interesting topics within enharmonicism are unaddressed; future research would seek to formalize concepts such as key signatures, modes, chords, triads, chord progressions, negative harmony, and chord substitutions. Of particular interest would be exploring how our enharmonic theory fits within the broader history of negative harmony, including Riemannian harmonic dualism, Levy's harmonic polarity, Coleman's symmetric interpretation of negative harmony, and Collier's reflection interpretation of negative harmony \cite{Haller}.

But even without addressing these additional topics, there are questions leading to interesting problems within this paper as it is. For example, what is the most efficient way to determine if an enharmonic system is RP-ME? To what extent are RP-ME enharmonic systems unique? Given a diatonic set, does there exist a formula for the elements of the $n$-tuple used to define enharmonic equivalence that guarantees the resulting enharmonic system is RP-ME? If such a formula were to exist, what properties would it need to satisfy, and if no such formula exists, could it be proven why not? Assuming it is possible to generate more than one RP-ME enharmonic system using distinct $n$-tuples and a fixed diatonic set, what would the set of such $n$-tuples look like, and what properties would it satisfy? These questions will be explored in a follow-up paper by the author.



\section*{Acknowledgements}
\addcontentsline{toc}{section}{Acknowledgements}
The author thanks the Lord Jesus Christ for his wisdom, guidance, grace, and mercy throughout the process of researching, writing, and revising. The author thanks those who took the time to review this article, offering generous feedback and detailed comments. The author also thanks Dr. Elizabeth Carlson for offering advice and encouragement during the research process and his parents Kyle and Kristine Alons for their unending love and support.

\newpage
\bibliographystyle{acm}
\bibliography{ref}

\addcontentsline{toc}{section}{References}

\end{document}